\documentclass[leqno]{article}
\usepackage{graphicx}
\usepackage{hyperref}
\usepackage{amsmath}
\usepackage{amssymb}
\usepackage{amscd}
\usepackage{mathrsfs}
\usepackage{comment}
\usepackage{dsfont}
\usepackage[T1]{fontenc}
\usepackage[utf8]{inputenc}
\usepackage{authblk}

\newtheorem{thm}{Theorem}

\newtheorem{lem}[thm]{Lemma}
\newtheorem{definition}{Definition}
\newtheorem{remark}{Remark}

\newcommand{\R}{\mathbb{R}}
\newcommand{\C}{\mathbb{C}}
\newcommand{\N}{\mathbb{N}}
\newcommand{\K}{\mathbb{K}}
\newcommand{\Z}{\mathbb{Z}}
\newcommand{\muz}{\mu}

\renewcommand{\u}{\psi}
\renewcommand{\v}{\varphi }

\newcommand{\scrap}[1]{}
\newcommand{\Zn}{\Phi}
\newcommand{\Lop}{\mathcal{L}}
\newcommand{\Mop}{\mathcal{M}}
\newcommand{\Nop}{\mathcal{N}}

\newcommand{\bydef}{:=} 

\renewcommand{\Re}{\textrm{Re}}
\renewcommand{\Im}{\textrm{Im}}
\newcommand{\scat}{\kappa}
\newcommand{\ar}{\sigma}
\newcommand{\cn}{\textrm{cn}}
\newcommand{\sn}{\textrm{sn}}

\def\eproof{{\hfill $\Box$ \vspace{12 pt}}}

\begin{document}

\title{Rigorous numerics for NLS: bound states, spectra, and controllability}

 \author[1]{Roberto Castelli\thanks{rcastelli@bcamath.org}}
\author[2,1]{Holger Teismann\thanks{hteisman@acadiau.ca}}
\affil[1]{\footnotesize{Basque Center for Applied Mathematics, Bilbao, Basque Country, Spain}}
\affil[2]{\footnotesize{Department of Mathematics and Statistics\\
  Acadia University, Wolfville, 
  Canada}}

 \date{\today}

\maketitle
\begin{abstract}
 In this paper it is demonstrated how rigorous numerics may be applied to the 
one-dimensional nonlinear Schr\"odinger equation (NLS); specifically, to determining bound--state solutions and 
establishing certain spectral properties of the linearization. Since the results are rigorous, 
they can be used to complete a recent analytical proof \cite{BLT12} of the local exact controllability of NLS.
\end{abstract}

{\small \textbf{Key words:} rigorous numerics,  radii polynomials, controllability of PDEs, spectral analysis, BEC
\\
 	}

\section{Introduction}
 Analytical proofs of interesting/important/desirable properties of mathematical models  
 are often asymptotic in nature; as such, they are liable to leave a  finite number 
 of cases
 undecided. Or an analytical argument may reduce the property of interest to a criterion that
 needs to be verified each time the property is to be established for a particular model. 
 In both cases, it is natural to resort to numerical methods to conclude the argument.If this reasoning is to have the stature of a rigorous proof, one needs to use rigorous numerics.\\ 
 A famous case in point is the proof of the existence of the Lorentz attractor by Tucker \cite{tucker1999lorenz,tucker2002rigorous};
 examples of this general scenario closer to the topic of the 
present paper include
\begin{itemize}
  \item the asymptotic formula for the ground-state energy of a non-relativistic atom
  (by numerical verification of an elementary inequality)
  \cite{feffermanSeco1996};
  \item the conditional asymptotic stability of solitary waves of the cubic nonlinear 
   Schr\"odinger equation (by numerical verification of the gap condition) \cite{schlag2006numerical};
   \item the absence of imbedded eigenvalues for linearized NLS \cite{simpson2011spectral,simpson2011embedded};
    \item the existence of surface gap solitons of the 1-D NLS 
  (by numerical verification of an integral inequality) \cite{plum2012interfaces};
  \item the exclosure of eigenvalues of the Schr\"odinger equation with a perturbed periodic potential
   \cite{plum2012eigenvalue}.
  \end{itemize}

In the present paper we apply this paradigm to spectral properties of the linearized NLS (on a finite 
interval with Dirichlet boundary conditions), which 
are needed in the analytical proof of its controllability \cite{BLT12}.  \\
 
 Controlled manipulation of quantum systems is a very active field in 
 science and engineering (several review papers and monographs are available on this 
 subject; for a recent survey, see e.g. \cite{brif2012control} and the literature therein). A class of systems that have been studied intensely are 
 Bose Einstein condensates (BEC).  In this paper we consider a one-dimensional condensate 
 in a hard-wall trap (``condensate-in-a-box"), where the trap size (box length) is a time-dependent 
 function $L(t)$ that can be manipulated. The precise definition of the model is 
 given in Section \ref{Def} below.
 Given an initial state, the objective is to ``engineer" the control function $L(t)$ such that
 the condensate will be guided to a particular target state.  
 The model we are considering was first proposed by Band et al. \cite{BMT02} to study adiabaticity in
 a nonlinear quantum system.  More recently, the opposite regime, fast transitions (``shortcuts to adiabaticity"),
  has been investigated for BECs in box potentials \cite{theodorakis2009oscillations,del2012shortcuts}.  
 Condensates in a box trap have actually been realized experimentally \cite{meyrath2005bose}, an accomplishment that 
 attracted considerable attention. 
 
In light of the these developments, it is natural to study the mathematical control properties of the nonlinear Schr\"odinger 
equation.  Ref. \cite{BLT12} establishes a local controllability result for  (\ref{NLSa}) in the vicinity of the nonlinear ground state
 $\phi $; the precise statement is given in Section \ref{Contr} below. 
      
  The proof relies on   two spectral properties of the operator $\Lop $ that arises by 
  linearizing equation (\ref{NLSA}) (the rescaled version of \eqref{NLSa}) about the state $\v(t,x) = e^{i\muz t} \phi (x)$; namely ($\Psi^{(1)}_n(x)$ denotes the first component 
  of the $n$-the eigenfunction of $\Lop ^*$;  see Section~\ref{Lin}), 
  \begin{enumerate}
 \item[(A)]  
  the integrals  $\Gamma _n = \int _0^1 (x\phi )'(x)  \overline{\Psi^{(1)}_n(x)} dx $
  are non-zero; 
   \item[(B)]   the non-zero eigenvalues $\lambda _n$ of $\Lop $ are simple. 
\end{enumerate} 
Analytical proofs of these properties are available, but the arguments are asymptotic in nature and are only applicable
for (potentially) large eigenvalue indices $n$. So for a finite number of cases the validity of (A) and (B)
is unclear (although it \textit{is} proved  that controllability holds \emph{generically};  see Section~\ref{Contr} for the precise statement).   By means of rigorous computation we close this gap; i.e. we give a rigorous 
computer-assisted proof that (A) and (B) hold \textit{for all} $n$.  Enclosure of    
 the nonlinear ground state and the spectrum of $\Lop $  is accomplished by 
applying \textit{radii polynomials} and suitable estimates for bi-- and trilinear convolution terms
arising from the nonlinearity. \\

This paper is organized as follows. 

In Section~\ref{Def} we state the original model 
and its rescaled version, which is the one we are working with throughout the paper. We 
discuss bound-state solutions (Section~\ref{Bound}) and describe the control problem and result 
proved in \cite{BLT12} (Section~\ref{Contr}). Section~\ref{NLSprob} closes with a description of the
linearization of NLS (around a given bound state) whose spectral properties will be studied by means of 
rigorous numerics. In Section~\ref{Comp} we give a general outline of the \textit{radii--polynomial method}  for performing rigorous
numerical computations.  How this general method is applied to the NLS-problem at hand is  described in Section~\ref{Appl}.   
Specifically, we describe the rigorous determination of bound states (Section~\ref{CompBound}) and eigenvalues and eigenfunctions 
of the linearization (\ref{CompEig}). We also explain our method for verifying that the eigenvalues  are simple (Section~\ref{CompSimple}). 
The final two sections of the main body of the paper contain an overview of the numerical results (Section~\ref{Num}) and 
some concluding remarks (Section~\ref{Dis}). There are two appendices at the end of the paper, which contain derivations of some 
the required estimates.

\section{Nonlinear Schr\"odinger equation and control problem} 
\label{NLSprob}
\subsection{Problem statement and rescaled equation}
\label{Def} 
Following Band et al. \cite{BMT02}, 
we consider the ``condensate-in-time-varying-box'' problem
\begin{subequations}
\begin{eqnarray}
 i\hbar \u _t &=& - \frac{\hbar ^2}{2m} \u _{xx}-\sigma\scat |\u |^2\u , \quad (x\in (0,L(t)),\, t\in (0,T)), \ \sigma=\pm 1
\label{NLSa} \\
 \u (t,0)&=& \u (t,L(t))=0 \quad (t\in [0,T]) \label{NLSb}
\label{NLS}
\end{eqnarray}
\end{subequations}
where 
\begin{enumerate}
 \item $\u (t,x )\in \C $ is the \textit{wave function}, which is assumed to be \textrm{normalized}; i.e.,
  \begin{equation}
   \| \u (t) \| _{\ell ^2(0,L(t))}^2 = \int _0^{L(t)} |\u (t,x)|^2 dx = 1.
   \label{norm}
   \end{equation} 
  In the control problem, $\u $ plays the role of the \textit{state}.
 \item $\hbar $ and $m$ are \textit{Planck's constant} and the \textit{particle mass};
 \item $\scat >0$  is a \textit{nonlinearity parameter}, derived from the scattering length and the particle number;
 \item the signs of  $\sigma$ correspond to the \textit{focussing} ($\sigma=1$) and   \textit{de-focussing} ($\sigma=-1$) cases, respectively; 
 \item $L_0>0$ is the initial (and final) length of the box (we will choose $L_0=1$ below); 
 \item $L:[0,T]\to (0,\infty )$ is a function such that $L(0)=L_0=L(T)$ and plays the role 
 of the \textit{control}. 
\end{enumerate}

\begin{remark}
i) This problem is a nonlinear variant of the control problem solved by 
K.~Beauchard \cite{B06}. 

ii) The normalization condition (\ref{norm}) can formally be derived from (\ref{NLSa}) \& (\ref{NLSb}), 
since the \textit{density}  $\rho (t,x) := |\u (t,x)|^2$ and the \textit{current} 
$J(t,x) := \frac{\hbar }{m} \Im (\bar{\psi } (t,x) \psi _x(t,x))$  satisfy the usual 
\textit{continuity equation} $\rho _t = -J_x $. 
\end{remark} 
To non-dimensionalize 
the problem and to transform it to the time-independent domain $(0,1)$, 
we introduce new variables \cite{BMT02}, 
\begin{equation}\u (t ,x ) := \frac{\hbar}{\sqrt{2\scat m} L(t)}\, \v \left( \frac{\hbar}{2m} \int _0^t \frac{ds}{L(s)^2} \, \,
\frac{x}{L(t)} \right) =: \frac{\hbar }{\sqrt{2\scat m} L(t)} \, \v (\tau ,\xi ). \label{trafo1}
\end{equation}

Moreover, defining $u(\tau ) = \frac{2m}{\hbar ^2} L(t) \dot{L} (t)$
and renaming $\tau \to t$, $\xi \to x$,  gives
\begin{subequations}
\begin{eqnarray}
 i\v _t &=& - \v _{xx}- \sigma|\v |^2\v -
 \textstyle i u(t)(x \v )_x  ,
 \quad (x\in (0,1),\, t \in (0,T))
\label{NLSA} \\
 \v (t,0)&=&\v (t,1)=0 \quad (t \in [0,T]). \label{NLSB}
\label{NLS*}
\end{eqnarray}
 \end{subequations} 
The spectral and control properties of (\ref{NLSA})-(\ref{NLSB}) are the subject of this paper.  
\subsection{Bound states}
\label{Bound}
We are looking for stationary solutions to 
the problem (\ref{NLSA}),(\ref{NLS*})  (with $u(\tau)\equiv 0$). To this end, let 
\begin{eqnarray*}
 \v (t ,x) &=& e^{i\ar\muz t}\phi (x)
\end{eqnarray*}
where\footnote{Note that the sign in the exponents of the time-dependent part of $\v $ depends on $\ar $. In the de-focusing case, the definition of $\muz $ is the 
one favoured by physicists. The choice of sign in the focusing case is 
consistent with \cite{RSS05}, which is one of our main references for analyzing the linearized equation.}  $\phi=\phi(x)$ is a
nonlinear bound state corresponding to the chemical potential  $\muz $; i.e. a real solutions of the boundary value problem
\begin{subequations}
\begin{eqnarray}
 -\phi'' &+&\ar \muz   \phi \,-\ar \phi^3=0 \quad (x \in (0,1))\label{groundA1} \\
 \phi(0) &=& \phi(1) =0. \label{groundB1} 
 \end{eqnarray}
\end{subequations}  

Explicit formulas for the 
solutions of (\ref{groundA1}),(\ref{groundB1})
are available in terms of Jacobian elliptic functions.
If $j\in \{0,1,2,\ldots \} $, then 
$\phi _j^\pm (x)$ will denote 
the (real-valued) solution of (\ref{groundA1}),(\ref{groundB1}) 
which possesses precisely $j$ zeros (``nodes'') within the interval
$(0,1)$. The node-less solution $\phi ^\pm:=\phi _0^\pm$ is referred to as the 
\textit{ground state};
the solutions $\phi _j^\pm $ ($j\ge 1$) with one or multiple nodes are called 
\textit{excited states.}  To find an explicit solution formula for $\phi _j^+$ and $\phi _j^-$,
we first solve the equation(s) 
\begin{subequations}
\begin{eqnarray}
 \textrm{focussing case ($\ar =1$)} \quad  \muz  &=& 4(j+1)^2(2k^2-1)K(k)^2, \,\,\, \muz  \in [-\mathring{\mu }_j ,\infty ) \label{k} \\
 \textrm{de--focussing case ($\ar =-1$)} \quad  \muz  &=& 4(j+1)^2(k^2+1)K(k)^2, \quad  \muz  \in [\,\mathring{\mu }_j ,\infty )  \label{l}
\end{eqnarray}  
\end{subequations} 
for $k$, where $K(k)$ denotes the complete elliptic integral 
of the first kind (see, e.g. \cite{AS65}) and $\mathring{\mu }_j :=(j+1)^2\pi ^2$. Note that,
since $K(k)$ is a strictly increasing continuous function of $k\in [0,1)$ 
satisfying $K(0) = \frac{\pi }{2} $ and $\lim _{k\to 1^-} K(k) =\infty $, equation (\ref{k}) [resp. (\ref{l})] has 
exactly one solution $k=k_j^+ (\mu )$ [resp.  $k=k_j^- (\mu )$] for any choice 
of parameters $\muz \in [-\ar \mathring{\mu }_j ,\infty )$  
and  $j\in \{ 0,1,2,\ldots \} $. Moreover, 
the functions  $k_j^\pm :  [-\ar \mathring{\mu }_j ,\infty )\to [0,1 )$ are continuous and strictly 
increasing as well, and satisfy $\lim _{s\to \infty} k_j^\pm (s) = 1$.
Writing $k_j^\pm =k_j^\pm (\mu  )$, the solutions $\phi _j^\pm $ of (\ref{groundA1}),(\ref{groundB1}) are given by 
\cite{CCR},\cite{CCR00}
\begin{subequations}
\begin{eqnarray}
  \phi _j^+(x) 
 &=& 2\sqrt{2} (j+1) k_j^+ K(k_j^+) \,\, \cn \Big(
 2(j+1) K(k_j^+) {\textstyle (x-\frac{1}{2} )} +[j]_2K(k_j^+),k_j^+\Big) ,
\label{phij}\\
\phi _j^-(x) &=& 2\sqrt{2}(j+1)k_j^-K(k_j^-)\,\,\sn \Big(2(j+1)K(k_j^-)x,k_j^-\Big)
\end{eqnarray}  
\end{subequations} 
where $\cn =\cn (x,k) $ and $\sn =\sn (x,k) $ are the Jacobian elliptic cosine and sine functions, respectively,  and $[j]_2 :=
j$ mod $2$. 

\subsection{Control problem} 
\label{Contr}
We now state the controllability result mentioned above, which roughly states that,
``generically'' w.r.t. the parameter $\muz $, 
exact controllability holds locally around the ground state. Here ``generic'' means the existence of an at most 
countable set  $J\subset (-\ar \pi ^2,\infty )$ of potentially exceptional $\muz $ values. 
Defining 
\[ \mathcal{H} := \{ f\in H^3(0,1;\C )\mid f(0) = f(1) = 0\} \quad \textrm{and} 
\quad \mathcal{S} := \mathcal{H} \cap \{ f\in L^2(0,1;\C )\mid \int_0^1 |f|^2dx =1\} ,\] 
the precise 
statement reads  
\begin{thm} \label{BLT} 
 \cite{BLT12} Let $\muz \in (-\ar \pi ^2,\infty )\setminus J$, $\phi = \phi _\muz $ the corresponding 
 ground state, and $T>0$. Then there exists a number $\delta =
 \delta (T,\muz ) >0$ such that 
 for all states $\v _0,\v_1\in \mathcal{S} $ satisfying 
 \[ \| \v _0 - \phi \| <\delta \quad \textrm{and} \quad  \| \v _1 - e^{i\ar \muz T}\phi \| <\delta \]
 there exists a control function $u\in H^1([0,T],\R )$ with 
 $\int _0^T u(t)dt =0 $ such that the unique solution 
 $\v \in C([0,T],\mathcal{H} )$ of (\ref{NLSA})-(\ref{NLSB}) satisfies $\v (0) =\v _0 $ and $\v (T) =\v _1 $. 
\end{thm}
The fact that the theorem cannot be stated for \emph{all} values of $\muz $ is due to the asymptotic nature of
the direct analytical proof of the properties (A) and (B), which only covers (potentially) large $n$.
Properties (A) and (B) can still be established generically (i.e. up to an at most countable set $J$ of 
possible exceptions) in an indirect way by using the analytic dependence of the operator $\Lop $ and 
its spectrum on the parameter $\muz$. However, while the genericity property implies that 
controllability holds with ``probability one w.r.t. \textit{random} choices" of $\muz $, for any \textit{particular}
 value of  $\muz $  the theorem can only be applied if (A) and (B) are verified for
  the finite number of cases  not covered by the direct  proof. 
  In the remainder of the paper we are going to demonstrate that this verification can be accomplished by rigorous numerical 
  computation. 
 
 \subsection{Linearization}
 \label{Lin}
 The proof of Theorem \ref{BLT} uses linearization around the ground state and the Implicit Function Theorem. 
If $\phi $ is a bound state, then the function $\varphi (t,x) = e^{i\ar \muz t}
\phi (x) $ is the unique solution of (\ref{NLSA})-(\ref{NLS*})
with $\varphi (0,x)=\phi (x) $. 
Now we linearize around $\varphi $. The result is    
\begin{subequations}
\begin{eqnarray}
 iz_t &=&-z_{xx} -\ar |\varphi |^2z -2\ar \Re (\varphi \bar{z})\varphi   \label{Va}\\
 z(t,0)&=&z(t,1)=0\ . \label{Vb} 
\end{eqnarray}
\end{subequations}
The time dependence of the term involving 
$\Re (...)$ is eliminated by the  
transformation  $e^{i\ar \muz  t} \tilde{z} (t) := z(t) $, 
which yields the BVP 
\begin{subequations}
\begin{eqnarray}
 i\tilde{z}_t &=&-\tilde{z}_{xx} +\ar \muz \tilde{z}-\ar\phi ^2\tilde{z}
 -2 \ar \phi ^2  \Re (\tilde{z})  \label{VVa}\\
 \tilde{z}(t,0)&=&\tilde{z}(t,1)=0 \label{VVb} 
\end{eqnarray}
\end{subequations}
It is natural to work with the real $(2\times 2)$-system arising from 
(\ref{VVa})-(\ref{VVb})
by decomposition in real and imaginary parts.
Consider the matrix operator 
\begin{equation}
 \Lop:= \left( \begin{array}{cc} 0 & -\Delta +\ar \muz  -\ar \phi ^2(x) \\
 \Delta -\ar \muz  +3\ar \phi ^2(x) & 0 \end{array} \right) 
 =:  \left( \begin{array}{cc} 0 & L_- \\
 -L_+& 0 \end{array} \right) ;
 \label{L}
 \end{equation}  
 ($\Delta $ denotes the one-dimensional ``Laplacian'' $\frac{d^2}{dx^2} $.)
 Then eq. (\ref{VVa}) takes the form
 \begin{equation}
  Z_t = \Lop Z  , \label{Z} 
 \end{equation} 
 where $Z(t,x) = {\Re(\tilde{z}(t,x)) \choose \Im (\tilde{z}(t,x))}$.
The operator $\Lop $ is the main object of study.
\subsubsection{Provable properties of the spectrum of $\Lop $ (and $\Lop ^*$) if $\phi $ is the
ground state}  
\begin{enumerate}
 \item The spectrum of $\Lop$ consists of eigenvalues only
 \item all non-zero eigenvalues $\{ \lambda _n, \bar{\lambda } _n\} _{n\ge 1}$
 are purely imaginary, i.e. 
\[ \lambda _n =  i\beta _n , \quad \bar{\lambda }_n= - i\beta _n,\quad \beta _n >0
\qquad (\forall n\ge 1) .\] 
\item the multiplicity of the eigenvalues is at most 2; 
 \item \label{simple} all, \textbf{but possibly finitely many},  non-zero eigenvalues are simple
 \item the multiplicity of the eigenvalue zero is 2; let 
  \[ \Zn_0^+ = {0\choose \phi } \quad \textrm{and} \quad 
 \Zn_0^- = {\partial _\mu \phi \choose 0} .\] 
Then 
  \begin{equation}
  \Lop \Zn_0^-=\Zn_0^+ \quad \textrm{and} \quad 
   \Lop \Zn _0^+ =0 . \label{jordan}
  \end{equation}
 and the vectors $\Zn_0^+$, $\Zn_0^-$ form a basis of the generalized null space for 
$\Lop$. 
  \item Notation: 
 $\Zn_1^+,\Zn_2^+,\ldots ,\Zn _1^-,\Zn_2^- ,\ldots $ denote the eigenvectors\footnote{Clearly, 
 these are unique up to normalization.}  
corresponding to the 
non-zero eigenvalues $\lambda _1,\lambda _2,\ldots $, 
$\bar{\lambda }_1,\bar{\lambda }_2,\ldots $; i.e., 
\begin{equation} 
\Lop \Zn _n^\pm = \pm i\beta _n\Zn_n^\pm  , \quad \beta _n>0 \qquad (n\ge 1)
\end{equation}   
and 
$\overline{\Zn _n^+} = \Zn _n^-$ for all $n\ge 1$ (where 
 $\overline{(.)}$ denotes complex conjugation).   
\item  Similarly, $\Psi _0^+,\Psi _0^-,\Psi _1^+,\Psi _2^+,\ldots ,{\Psi }_1^-,{\Psi }_2^-,\ldots $ denote the   
  eigenfunctions for $\Lop ^*$ with corresponding eigenvalues 
$\bar{\lambda} _n$ and ${\lambda }_n$,
  respectively; i.e., 
  \begin{equation}
  \Lop ^* \Psi_n^\pm = \mp i\beta _n\Psi_n^\pm ,\quad \beta_n >0 \qquad
   (n\ge 1 ) 
   \end{equation}  
   Moreover, 
  \[ \Lop ^*\Psi _0^+ = \Psi _0^- \quad \textrm{and} \quad
   \Lop ^* \Psi _0^- =0 .\] 
where  
\begin{subequations}
 \begin{eqnarray}
 \Psi _0^- &=& {\phi \choose 0 } \label{Pa} \\
 \Psi _0^+ &=& { 0 \choose \partial _\mu \phi } \label{Pb}\\
\overline{\Psi_n^+} &=& \Psi _n^-
 \end{eqnarray}
 \end{subequations} 
 \item $\{ \Zn _m^\pm \}_{m\ge 0 }$, $\{ \Psi _n^\pm \}_{n\ge 0 }$
 form  bi--orthogonal systems; i.e., 
  \begin{eqnarray}
  \langle \Zn_m^\sigma ,\Psi _n^\tau \rangle   &=& \delta _{m,n}^{\sigma , \tau} , \qquad 
  m,n\in \{ 0,1,2,\ldots \} \label{biA} ,\,\, \sigma ,\tau \in \{ +,-\}   
 \end{eqnarray}
where the inner product $\langle.,.\rangle $ is defined by
\begin{eqnarray} 
 \langle U,V\rangle &=& \left\langle  \left( \begin{array}{c} U^{(1)} \\ U^{(2)} 
 \end{array} \right) ,\left( \begin{array}{c} V^{(1)} \\ V^{(2)} 
 \end{array} \right) \right\rangle    
  =   \int _0^1U^{(1)}(x)\overline{V^{(1)}(x)} dx
 +  \int _0^{1}U^{(2)}(x)\overline{V^{(2)}(x)} dx . \nonumber
\label{bracket}
\end{eqnarray}
and 
\[ \delta _{m,n}^{\sigma , \tau} = \left\{ \begin{array}{cl} 1, & m=n \, \, \textrm{and} \, \, \sigma = \tau \\
 0, & \textrm{otherwise} \end{array}. \right. \] 
\end{enumerate}
\begin{remark} Note that (other than in Section~\ref{Bound} above) the $\pm $ superscripts do \textit{not} refer to the focussing and defocusing cases here.  
Note also that the eigenfunctions  $\Zn_1^+,\Zn_2^+,\ldots ,{\Zn} _1^-,{\Zn}_2^- ,\ldots $ and
 $\Psi _1^+,\Psi _2^+,\ldots ,{\Psi }_1^-,{\Psi }_2^-,\ldots $ are complex-valued. 
 \end{remark}
\subsubsection{A change of variables} 
It is convenient to employ a similarity transformation
\cite[(12.15)]{RSS05}: Let
\[ J:= \left( \begin{array}{cr} 1 & i\\
 1 & -i \end{array} \right)  .\]
Then 
\begin{subequations}
\begin{eqnarray*} i\Lop &=& J^{-1}
\mathcal{M} J \label{J} , \quad 
 -i\Lop^* = J^{-1}
 \mathcal{N} J \qquad \textrm{and so}
 \end{eqnarray*} 
 \end{subequations}
\[ \textrm{spec}(\Lop) = i\,\textrm{spec}(\Mop) ,\quad  \textrm{spec}(\Lop ^*) = -i\,\textrm{spec}(\Nop) , \]  
where
\begin{eqnarray*}
\Mop &:= & \left( \begin{array}{rr} -\Delta   &  0 \\
 0 & \Delta     \end{array} \right) +  \sigma\left( \begin{array}{cc} \mu
 -2\phi ^2& -\phi ^2\\
 \phi ^2& -\mu +2\phi^2   \end{array} \right)   \\
\Nop &:= & \left( \begin{array}{rr} -\Delta   &  0 \\
 0 & \Delta     \end{array} \right) +  \sigma\left( \begin{array}{cc} \mu
 -2\phi ^2& \phi ^2\\
 -\phi ^2& -\mu +2\phi^2   \end{array} \right). 
 \end{eqnarray*}   
Now let  $(\pm \beta _n,V_n^\pm )$ and $(\pm \beta _n,W_n^\pm )$, $\beta _n\ge 0$,  
be the eigenpairs for the operators $\Mop $ and $\Nop $, 
respectively,  i.e., for  $n\ge 1$, $\beta _n >0$, 
\begin{subequations}
\begin{eqnarray}
 \Phi_n^\pm &=& J^{-1}V_n^\mp, 
 \quad \quad \,\,\,
\Psi_n^\pm = J^{-1}W_n^\mp,  
\label{M2} \\
  \Mop V_n^\pm &=& \pm \beta _nV_n^\pm , \quad 
  \Nop W_n^\pm \, = \pm \beta _nW_n^\pm . \label{M3}
\end{eqnarray}  
\end{subequations}
Writing $\beta =\pm \beta _n$, $V=V_n^\pm ={u\choose v}$, 
 $W=W_n^\pm ={w\choose z}$,  the characteristic 
equations (\ref{M3}) are equivalent to the BVP
\begin{subequations}
\begin{eqnarray}
 u'' -(\sigma\muz  -\beta ) u &=& -\sigma\phi ^2 (2u+v), \qquad \,u(0)=u(1)=0 \label{a}\\
 v'' -(\sigma\muz  +\beta ) v &=& -\sigma\phi ^2 (u+2v), \qquad v(0)=v(1) =0 \label{b} \\
 w'' -(\sigma\muz  -\beta ) w &=& -\sigma\phi ^2 (2w-z), \qquad \,w(0)=w(1)=0 \label{c}\\
 z'' -(\sigma\muz  +\beta ) z &=&  -\sigma\phi ^2 (-w+2z), \qquad z(0)=z(1) =0 \label{d} .
\end{eqnarray}   
\end{subequations}

\section{Computational method: rigorous computation using radii polynomials}
\label{Comp}
In this section we describe the rigorous computational method that will be used to 
\begin{itemize}
\item [(P1)] enclose the function $\phi(x)$ solution of \eqref{groundA1},\eqref{groundB1}; 
\item [(P2)] enclose the eigenpairs $(\beta, W)$ solutions of \eqref{c},\eqref{d};
\item [(P3)] prove that the eigenvalues $\beta$ are simple. 
\end{itemize}
These computations are based on suitable adaptations of the general method  known as {\it radii polynomials}. The radii--polynomial approach,  first introduced in \cite{MR2338393}, aims at demonstrating existence and local uniqueness of solutions of nonlinear problems by  verifying the hypothesis of the contraction mapping theorem in Banach spaces. In recent years this technique has been successfully applied to a variety of nonlinear problems; see e.g \cite{MR2443030, MR2821596, floquet, enclosure} and the references therein.

Before outlining the main steps of the method, as applied to P1, P2, and P3 above, we introduce some notation. For $z\in \C$, denote by $|z|=\max\{|\Re(z)|, |\Im(z)|\}$ and for a matrix $V=\{V_{i,j}\}\in \C^{n\times m}$ denote by $|V|=\{|V_{i,j}|\}$ and $|V|_{\infty}=\max_{i,j}|V_{i,j}|$. Given two matrices $A,B\in \R^{n\times m}$ the inequality $A<B$ is to be interpreted componentwise, i.e. $A_{i,j}<B_{i,j}$, for all $i,j$.
Define the weights $w_{k}$ as
$$
 w_{k}=\left\{\begin{array}{ll}
 1,\quad &k=0\\
 |k|,\quad &k\neq 0
 \end{array}\right.\ .
 $$
 
If $x=\{x_{k}\}_{k\geq 0}\in (\K^{d})^{\N}$ is a  sequence in $\K^{d}$ (for some $d\geq 1$ and $\K=\C$ or $\K=\R$) and $s>0$, we define the $s$-norm of $x$ as
 $$
\|x\|_{s}=\sup_{k}\{|x_{k}|_{\infty}w_{k}^{s}\}
$$ 
and $X^{s}$ the space 
$$
X^{s}=\{x\in (\K^{d})^{\N}\ |\ \|x\|_{s}<\infty \}.
$$ 
The space $(X^{s},\|\cdot\|_{s})$ is a Banach space; 
we refer to $s$ as the {\it decay rate parameter}. Let $B(r)=\{ x\in X^{s}\ |\ \|x\|_{s}\leq r\}$ be the ball of radius $r$ in $X^{s}$ and, for any $x\in X^{s}$, denote by 
\begin{equation}\label{balls}
B_{x}(r)=x+B(r)
\end{equation} 
the ball centred at $x$.

The first step of the method consists of rephrasing the original problem in terms of  an equation of the form 
$$
f(x)=0,
$$
where $f:X\to W$ is a  (possibly) nonlinear operator with  $X=X^{s_{1}}$, $W=X^{s_{2}}$ and suitable $s_{1},s_{2}$.
 
 Then we choose the {\it finite dimensional parameter} $m\ge 1$ and define the finite dimensional projections
$ \Pi_m^X : X \rightarrow X_m$ and $   \Pi_m^W : W \rightarrow W_m$ 
 as well as the infinite ``tail projections"
$ \Pi_\infty^X : X \rightarrow X_\infty$ and $\Pi_\infty^W : W \rightarrow W_\infty$ 
by
$$
\Pi_m^X(x)=x^{(m)}=(x_{0},\dots, x_{m}),\quad \Pi_\infty^X(x)=x^{\infty}=(x_{m+1},x_{m+2},\dots),
$$
and similarly for $   \Pi_m^W$ and $\Pi_\infty^{W}$. 

Consider the finite dimensional projection of the map $f$,
\begin{equation}
\begin{split}
f^{(m)}:&X_m \rightarrow W_m\\
& x \mapsto f^{(m)}(x) \bydef \Pi_m^W f(x,0^{\infty}),
\end{split}
\end{equation}
and suppose that an approximate solution $\bar x\in X^{(m)}$ of $f^{(m)}(x)=0$ has been computed. Sightly abusing notation,  
we use  $\bar x $ to indicate both the vector in $X^{(m)}$ and the sequence $(\bar x,0^{\infty})\in X$. Hence we also refer to $\bar x$ as the approximate zero of the full (infinite dimensional) map $f$, i.e.
$$
f(\bar x)\approx 0 .
$$
The next step is to define a nonlinear operator $T:X\to X$ with the property that the zeros of $f(x)$ are in one--to--one correspondence with the fixed points of $T$. The fixed point operator $T$ will be defined as a modified  Newton operator centred at the numerical  solution $\bar x$. The crux of the method is to prove that $T$ is a contraction.

Let 
$$
Df^{(m)}\bydef\frac{\partial f^{(m)}}{\partial x^{(m)}}(\bar x)
$$
be the Jacobian of $f^{(m)}$ evaluated at $\bar x$,  
$$
\Lambda_{k}=\frac{\partial f_{k}}{\partial x_{k}}(\bar x) \qquad (k=0,1,2,\ldots ),
$$
and   $A^{(m)}\in \K ^{(m+1)\times (m+1)}$ an invertible  approximate inverse of $Df^{(m)}$. Then we define the operator $A$ by  
\begin{equation}\label{eq:A}
(Ax)_{k}\bydef\left\{\begin{array}{ll}
(A^{(m)}x^{(m)})_{k}, &  k\leq m\\
\Lambda_{k}^{-1}x_{k},&k> m
\end{array}\right.
\end{equation}
and the fixed point operator $T:X\to X$ as
\begin{equation}\label{eq:T}
T(x)=x-Af(x).
\end{equation}
To ensure that fixed points for $T$ correspond to zeros of $f(x)$, we need to prove that the operator $A$ is injective. 
Since the finite part $A^{(m)}$ is invertible by construction, this amounts to verifying  that for $k>m$ the operators $\Lambda_{k}$ are invertible as well. 

The existence (and uniqueness) of the fixed point for $T$ will follow from Banach's fixed point theorem, once the operator $T$ has been proven to be a contraction on a suitable subset of $X$. The candidate sets are the balls  $B_{\bar x}(r)$ defined in \eqref{balls}, hence we  need rigorous estimates for the image of $T$ and  the rate of contractivity of $T$ on these balls.
 
Suppose that, for a fixed \textit{computational parameter} $M$, we have found bounds $Y=\{ Y_{k}\}_{k<M}$, $Z=\{Z_{k}(r)\}_{k<M}$, $Y_{M}$, and $Z_{M}(r)$, such that
\begin{equation}\label{eq:boundsYZ}
|(T(\bar x)-\bar x)_k|\leq Y_{k},\quad\sup_{b_1,b_2\in B(r)}\Big|\big[DT(\bar x+b_1)b_2\big]_{k}\Big|\leq Z_{k}(r)\qquad (\forall k<M), 
\end{equation}
and 
\begin{equation}\label{eq:codaM}
|(T(\bar x)-\bar x)_k|_{\infty}\leq \frac{ 1}{w_k^{s}}Y_{ M},\quad\sup_{b_1,b_2\in B(r)}\Big|\big[DT(\bar x+b_1)b_2\big]_{k}\Big|_{\infty}\leq\frac{1}{w_k^{s}} Z_{M}(r)\qquad (\forall k\geq M).
\end{equation}
\begin{definition} \label{def:radii_polynomials}
The {\em radii polynomials} are defined as
\begin{equation}
\begin{array}{l}
p_{k}(r)\bydef |Y_{k}+Z_{k}(r)|_{\infty}-\frac{r}{{w_k}^{s}} ,\quad k=0,\dots, M-1,\\
p_{ M}(r)\bydef Y_{ M}+ Z_{ M}(r)-r . \label{radpoly}
\end{array}
\end{equation}
\end{definition}
These are called {``polynomials"} because each bound $Z_{k}$ will be constructed  as a polynomial in the variable $r$ with degree equal to the degree of nonlinearity of the map $f(x)$. 

The final step in the procedure is  to solve the inequalities $p_{k}(r)<0$ for $r$. Then $T$ will be a contraction in any ball $B_{\bar x}(r^{*})$ whose radius $r^*$ satisfies $p_{k}(r^{*})<0$ for all $k\in \{ 0,1,\ldots M\} $. This is the content of the next
theorem. 

\begin{thm} \label{thm:rad_poly}
Suppose  $Y=\{Y_{k}\}_{k}, Z(r)=\{Z_{k}(r)\}_{k}$  satisfy \eqref{eq:boundsYZ} for $k=1,\dots, M-1$ and $Y_{M}$, $Z_{M}$ satisfy \eqref{eq:codaM} and let 
the polynomials $p_{k}(r)$, $p_{M}(r)$ be defined by \eqref{radpoly}.
Then, for every number $r>0$ such that $p_{k}(r)<0$ for all $k=0,\dots, M $, there exists a unique $x^* \in B_{\bar x}(r)$ such that $f(x^*)=0$.
\end{thm}
{\it Proof}.
See \cite{MR2338393}. \hfill $\Box $ \\

The enclosure radius $r$ arises as a solution of $p_{k}(r)<0$ ($k=0,\dots, M $), where the polynomials $p_{k}(r)$ are constructed from 
analytical estimates and numerical computations. Although the method relies on computer calculations, the results are mathematically rigorous because  all computations are performed in interval arithmetic (using the software package INTLAB \cite{Ru99a}), which accounts for all possible rounding errors.   

In addition to the radii--polynomial technique, there are several other computational methods based on the Contraction Mapping Principle (CMP), such as the    Krawczyk operator approach \cite{MR0255046, MR2394226}  or the methods developed by Yamamoto \cite{MR1639986}, by Koch et al. \cite{MR2679365}, and by Nagatou et al. \cite{Nagatou2011uq}. The main difference is that in the the radii--polynomials approach the enclosure radius $r$ is computed a posteriori and optimally, whereas in most of the other methods an initial guess is made of the set on which $T$ might be contractive, and the hypotheses of the CMP are verified after the fact. We feel that our approach has at least two advantages: the use of interval arithmetic is deferred to the end of the process reducing computing time, and the procedure attempts to determine an enclosure radius
that is as small as possible. The second consideration is particularly relevant to this work: the computation of the spectrum of $\mathcal{L} $ requires prior computation of $\phi(x)$ and the size of the intervals can growth dramatically when a large number of interval computations is performed; it is therefore necessary to have a very narrow enclosure of the solution $\phi(x)$. 
To accomplish this, we need sharp analytical estimates 
to control the truncation error arising from the finite dimensional approximation. 
 
In  summary, the technique consists of the following steps:
\begin{itemize}
\item[1.] to formulate the problem in the form $f(x)=0$ for a suitable map $f:X\to W$;
\item[2.] to fix a finite--dimensional projection, compute a numerical solution $\bar x$, and construct the fixed point operator $T$;
\item[3.] to compute the bounds $Y_{k}$, $Z_{k}$, $Y_{M}$ and $Z_{M}$ and construct the radii-polynomials;
\item[4.] to determine $r$ such that $p_{k}(r)<0$.
\end{itemize}

\subsection{Construction of the radii polynomials}\label{sec:radiipolynomials}
 
The construction of the bounds $Y$ and $Z$ is described next. 
 First we fix a  computational parameter $M$, ($M>m$), and we compute a constant $C_{\Lambda}$ so that 
 \begin{equation}\label{eq:CLambda}
\|\Lambda_{k}^{-1}\|_{\infty}\leq C_{\Lambda}\quad \forall k\geq M\ .
\end{equation}

Since $T(\bar x)-\bar x=Af(\bar x)$,  we define
\begin{equation}\label{eq:Y}
Y_{k}\bydef\left\{\begin{array}{ll}
|[A^{(m)}f^{m}(\bar x)]_{k}|,\quad &k\leq m\\
|\Lambda_{k}^{-1}f_{k}(\bar x)| &m+1\leq k\leq M-1\\
\end{array}\right. .
\end{equation}
In order to construct  the bound $Z_{k}$, we introduce the  operator 
$$
(J^{\dag}x)_{k}\bydef\left\{\begin{array}{ll}
(Df^{(m)}x^{m})_{k}, &  k\leq m\\
\Lambda_{k}x_{k},&k> m
\end{array}\right.
$$
and consider the splitting
\begin{equation}\label{eq:split}
\begin{array}{rl}
DT(\bar x+b_{1})b_{2}=&\left[I-ADf(\bar x+b_{1}) \right]b_{2}\\
=&\left[I-AJ^{\dag}\right]b_{2} - A\left[Df(\bar x+b_{1})- J^{\dag}\right]b_{2}
\end{array}\ .
\end{equation}
Since $b_{1},b_{2}\in B(r)$, it is convenient to write $b_{1}=ru$, $b_{2}=rv$, with $u,v\in B(1)$ and from the previous formula we have
\begin{equation}\label{eq:splitnorm}
\Big|\left[DT(\bar x+ru)rv\right]_{k}\Big|\leq_{cw}\Big|\left[\left(I-A J^{\dag}\right)rv\right]_{k}\Big| +\Big|\left[ A \left(Df(\bar x+ru)- J^{\dag}\right)rv\right]_{k}\Big|.
\end{equation}
Let $Z^{0}$ be defined as  
\begin{equation}
 (Z^{0})_{k}=\left\{\begin{array}{ll}
\left[|I-A^{m}Df^{(m)}|\{w_{j}^{-s}\}_{j\leq m}\right]_{k}, \quad &k\leq m\\
~ 0, & k> m
\end{array}\right.\ 
\end{equation} 
so that $\Big|\left[\left(I-A J^{\dag}\right)rv\right]_{k}\Big|\leq Z^{0}_{k}r$.

According with the degree $p$ of nonlinearity of the function $f(x)$, we can expand  $\left[(Df(\bar x+ru)- J^{\dag})rv\right]_{k}$ as a polynomial in $r$ 
\begin{equation}\label{eq:expansion}
\left[\left(Df(\bar x+ru)- J^{\dag}\right)rv\right]_{k}=\sum_{i=1,\dots,p}c_{k,i}r^{i}
\end{equation}
and we define the bounds $Z_{k}^{i}$ so that $Z_{k}^{i}\geq |c_{k,i}|$ uniformly in $u,v\in B(1)$. Finally the bound $Z_{k}$ is given by
\begin{equation}
 Z_{k}\bydef\left\{\begin{array}{ll}
[|A^{m}|((Z^{1})^{m}r+(Z^{2})^{m}r^{2}+\dots+(Z^{p})^{m}r^{p}]_{k}+Z^{0}_{k}r& k\leq m\\
|\Lambda_{k}^{-1}|(Z^{1}_{k}r+Z^{2}_{k}r^{2}\dots+Z^{p}_{k}r^{p})& m+1\leq k<M\\
\end{array}\right.\ .
\end{equation} 
Here $(Z^{i})^{m}=\Pi_{m}^{X}(Z^{i})$, which is the vector with the components $Z^{i}_{k}$ for $k\leq m$. 

The definition of the tail bounds $Y_{M}$ and $Z_{M}$ satisfying \eqref{eq:codaM} follows from uniform estimates, up to $w_{k}^{-s}$, of $|f_{k}(\bar x)|$ and $|c_{k,i}|$ for $k\geq M$, where we assume to have found $f_{M}$, $Z_{M}^{i}$ such that
\begin{equation}\label{eq:unifest}
|f_{k}(\bar x)|_{\infty}\leq \frac{1}{w_{k}^{s}}f_{M},\qquad |c_{k,i}|_{\infty}\leq \frac{1}{w_{k}^{s}}Z_{M}^{i} ,\qquad \forall k\geq M,\quad  \forall i=1,\dots,p .
\end{equation}
Then, in view of \eqref{eq:CLambda}, we define 
$$
Y_{M}\bydef C_{\Lambda}f_{M} \qquad Z_{M}\bydef C_{\Lambda}(Z^{1}_{M}r+\dots+Z^{p}_{M}r^{p}).
$$
We remark that the definition of the vector $Y$ and $Z$ is based on a  combination of rigorous computations and analytical estimates: utilizing rigorous computation ensures that the rounding errors are controlled whenever a computation is performed;
analytical estimates control the truncation errors arising from the finite--dimensional approximation. In particular, analytical estimates will be necessary to control 
 $f_{k}(x)$ for $k\geq M$ and to bound  the coefficients $c_{k,i}$ appearing in \eqref{eq:expansion}, both for each $k<M$ and uniformly for $k\geq M$.

\section{Application to NLS}
\label{Appl}
We now apply the computational technique described in the previous section to the control problem of Section~\ref{Def}. As mentioned in the introduction, the goal is to prove conditions (A) and (B) for \textit{all} eigenvalue of the linearized NLS.  To check whether the $\Gamma_{n}$ are non-zero, we first have to rigorously compute the eigenvalues and the eigenfunction of $\mathcal L$, i.e. the  solutions of system \eqref{c},\eqref{d}.  Since the linearization depends on the solution $\phi(x)$ of the Schr\"odinger equation, the bound state $\phi(x)$ has to be rigorously computed as well. 
Hence we perform three computations, each one using rigorous numerics.  
\begin{itemize}
\item[i)] For a choice of $\mu>0$ and $\sigma\in \{\pm 1\}$, we compute the solution $\phi(x)$ for  \eqref{groundA1}, \eqref{groundB1};
\item[ii)] Given the state $\phi(x)$, we compute the eigenpairs   $(\beta,W)$ by solving  \eqref{c}, \eqref{d} and check that $\Gamma_n$ is different from zero;
\item[iii)] We prove that the computed eigenvalues are simple.
 \end{itemize}
 
For each of these problems we state the nonlinear map $f(x)$, the Banach space $X^{s}$, the Jacobian $Df^{(m)}$, and some of the necessary analytic estimates. However, in order to increase readability, we delegate most of the analytical estimates and the technical details  to the Appendix.

\subsection{Computing the bound states $\phi(x)$}
\label{CompBound}
Bound--states $\phi$ are solutions of the BVP
\begin{equation}\label{eq:BVP}
\left\{\begin{array}{l}
-\phi''+\sigma\mu\phi-\sigma\phi^{3}=0,\quad x\in(0,1)\\
\phi(0)=\phi(1)=0
\end{array}\right.\ .
\end{equation}
Expanding $\phi$ w.r.t the sine-basis  $\{\sqrt 2 \sin(\pi n x)\}_{n\geq 1}$ gives 
$$
\phi(x)=\sqrt 2\sum_{n\geq 1}\alpha_{n}\sin(\pi n x),\qquad \alpha_{n}\in \R.
$$
Using the symmetry of the sine functions, this expansion is equivalent to 
\begin{equation}\label{eq:phiinbetas}
\phi(x)=\sqrt 2\sum_{n\in\Z}b_{n}\sin(\pi n x),\qquad b_{n}\in \R ,
\end{equation}
where the coefficients satisfy $ b_{-n}=-b_{n}$. This is readily seen by defining $\alpha_{n}=2b_{n}$.
 The advantage of this representation is that the projection of the cubic term onto the basis elements $\sqrt 2\sin(\pi n x)$ takes the  simple form
\begin{equation}\label{eq:cubica}
<\phi^{3},\sqrt 2\sin(\pi n \bullet)>\,=-4\sum_{\substack{
p+k+\ell=n \\
p,k,\ell\in \Z
}}b_{p}b_{k}b_{\ell};
\end{equation}
see Appendix A.
Inserting  \eqref{eq:phiinbetas}, \eqref{eq:cubica} into  system \eqref{eq:BVP} 
and using $b_{-n}=-b_{n}$,   
we obtain the infinite--dimensional algebraic system 
$$
f(b)=(f_{1},f_{2},\dots,f_{n},\dots)(b)=0, \quad n\geq 1 
$$
for the unknown $b=\{b_{n}\}_{n\geq 0}$, where

\begin{equation}\label{eq:gn}
f_{n}(b)=(\pi^{2}n^{2}+\sigma\mu)b_{n}+2\sigma\sum_{\substack{
p+k+\ell=n \\
p,k,\ell\in \Z
}}b_{p}b_{k}b_{\ell}\qquad n\geq 1\ .
\end{equation}
Note that we only considered $n\geq 1$: by the symmetry of the $b_{k}$'s we have that $f_{-n}(b)=-f_n(b)$. Since the unknowns are $b_{k}$ with $k\geq 1$ ($b_{0}$ may be set equal to zero), it is sufficient to solve $f_{n}(b)=0$ for $n\geq 1$.

\subsubsection{Ground State and even exited states}

The ground state and the exited states with an even number of nodal points are functions that are symmetric with respect to $x=\frac{1}{2}$. This means that the even Fourier coefficients vanish, i.e. $\alpha_{2n}=0$ and $b_{2n}=0$, as well as $f_{2n}(b)=0$.
Hence we discard the even Fourier coefficients and we introduce the sequence of {\it odd} coefficients
$$
b^{o}_{n}=b_{2n-1}.
$$
The symmetry conditions for the new sequence read  $b^{o}_{0}=-b^{o}_{1}$ and $b^{o}_{-n}=-b^{o}_{n+1}$ for $n\geq 1$.

Similarly, we discard the even component of $\{f_{n}\}$ and introduce the reduced system $f^{o}_{n}=f_{2n-1}$ for $n\geq 1$. 
In  terms of the unknown  $b^{o}=\{b^{o}_{n}\}_{n\geq 1}$ the new system  reads 
$$
f^{o}_{n}(b^{o})=(\pi^{2}(2n-1)^{2}+\sigma\mu)b^{o}_{n}+2\sigma\sum_{\substack{
p+k+\ell =n+1 \\
p,k,\ell\in \Z
}}b^{o}_{p}b^{o}_{k}b^{o}_{\ell } .
$$
We wish to bound the solution $b^{o}=\{b^{o}_{k}\}_{k\geq 1}$ of $f^{o}(b^{o})=0$. For the remainder of this section we omit the superscript $(\cdot) ^{o}$.
We look for the solution in the Banach space
$$
X^{s}=\{b=\{b_{k}\}_{k\geq 1}, b_{k}\in \R : \|b\|_{s}<\infty \}
$$
for $s\geq 2$. Note that $f:X^{s}\to X^{s-2}$.
  
Suppose that the finite--dimensional parameter $m$ has been chosen and that a numerical solution $\bar b=\{ \bar b_{1},\bar b_{2},\dots, \bar b_{m}\}$ of $f^{(m)}(b)=0$ has been computed. (A package such as  {\it Maple} may conveniently be used  to determine the Fourier coefficients of the elliptic functions up to a desired accuracy.)

By direct computation, the Jacobian of $f^{(m)}$  and the coefficients $\Lambda _n$ are  given by 
\begin{equation}
\begin{split}
\frac{\partial f_{n}}{\partial b_{j}}(\bar b)=(\pi^{2}(2n-1)^{2}&+\sigma\mu)\delta_{n-j}+6\sigma\Big[\sum_{k_{1}+k_{2}=n-j+1}\bar b_{k_{1}}\bar b_{k_{2}}-\sum_{k_{1}+k_{2}=n+j}\bar b_{k_{1}}\bar b_{k_{2}}\Big],\quad  k_{1},k_{2}\in\Z ,
\end{split}
\end{equation}
and 
$$
\Lambda_{n}\bydef\pi^{2}(2n-1)^{2}+\sigma\mu, \quad  n>m,
$$
respectively. 
We now introduce the operator $T$ according to   \eqref{eq:A}, \eqref{eq:T}. (In the de-focusing case $\sigma=-1$ the parameter $m$ must be such that   $\pi^{2}(2m-1)^{2}>\sigma\mu$ to ensure the invertibility of $\Lambda_{n}$ and, by extension, of $A$). 
Note that the Jacobian is symmetric, as expected from the variational nature of the problem. 

For the  construction of the radii polynomials we fix the computational parameter $M=3m$ and set
$$
C_{\Lambda}=(\pi^{2}(2M-1)^{2}+\sigma\mu)^{-1}.
$$
The choice of $M$ is motivated by the fact that 
 $f_{n}(\bar b)=0$ ($\forall n\geq 3m$), which allows us to set $Y_{M}\bydef 0$.  The definition of the vector $Y=(Y_{1},\dots,Y_{M-1})$ is given by \eqref{eq:Y}, while the vector $Z=(Z_{1},\dots, Z_{M-1})$ and the tail bound $Z_{M}$ follow from careful estimates  of  the coefficients $c_{k,i}$, given in terms of convolution products. We use the estimates provided in the paper \cite{MR2718657}, where 
 sharp bounds for the convolution products  are  proved. In Appendix A we list some of required analytical estimates and the definition of the remaining bounds $Z_{k}, Z_{M}$; see \cite{MR2718657} for details.

\subsubsection{Odd exited states}

The procedure for  the computation of the odd exited states (i.e. solutions $\phi(x)$ with an odd number of nodal points)  is similar to the one discussed in the previous section. Since the solutions are odd w.r.t. $x=\frac{1}{2}$,  only the \textit{even}  Fourier coefficients have to be computed.  Consequently,  we introduce the vector of unknowns $b^{e}_{n}=b_{2n}$ and the system
$$
f^{e}_{n}(b^{e})=(\pi^{2}(2n)^{2}+\sigma\mu)b^{e}_{n}+2\sigma\sum_{\substack{
p+k+\ell =n \\
p,k,\ell\in \Z
}}b^{e}_{p}b^{e}_{k}b^{e}_{\ell }
$$
to be solved for $n\geq 1$. 
Note that in this case the Jacobian is given by
\begin{equation}
\begin{split}
\frac{\partial f^{e}_{n}}{\partial b^{e}_{j}}(\bar b)=(\pi^{2}(2n)^{2}&+\mu)\delta_{n-j}+6\Big[\sum_{k_{1}+k_{2}=n-j}\bar b_{k_{1}}\bar b_{k_{2}}-\sum_{k_{1}+k_{2}=n+j}\bar b_{k_{1}}\bar b_{k_{2}}\Big],\quad  k_{1},k_{2}\in\Z.
\end{split}
\end{equation}
The definition of the fixed point  operator $T$ and the construction of the radii polynomials are  similar, mutatis mutandis, to the previous case.

\begin{remark}\label{rmk:enclosure}
For clarity, we explicitly show what the enclosure of the sequence $b^{o}$ or $b^{e}$ means for the actual Fourier coefficients $b_{n}$ in  \eqref{eq:phiinbetas}. Depending on the symmetry of the state $\phi(x)$, denote by $b^{*}_{n}$ the odd $b^{o}_{n}$  or the even  $b^{e}_{n}$  coefficients.

Suppose that, for a finite--dimensional parameter $m=m_{\phi}$ and a decay rate $s=s_{\phi}>2$, the computational method results in the enclosure radius $r=r_{\phi}$.  This means that the sequence  $b^{*}=\{b^{*}_{n}\}_{n\geq 1}$ satisfies
$$
|b^{*}_{n}-\bar b^{*}_{n}|\leq r_{\phi}/w_{n}^{s_{\phi}} \ {\rm for }\  n=1,\dots, m_{\phi},\quad {\rm and}  \quad 
|b^{*}_{n}|\leq r_{\phi}/w_{n}^{s_{\phi}} \ {\rm for }\  n>m_{\phi} .
$$
Hence, the sequence $b=\{b_{n}\}$ satisfies
\begin{equation}\label{eq:enclosurebn}
|b_{n}-\bar b_{n}|\leq r_{\phi}\frac{1}{w_{\left[\tfrac{n}{2}\right]}^{s_{\phi}}} \ {\rm for }\  |n|\leq 2m_{\phi}, \quad 
 {\rm and}  \quad |b_{n}|\leq r_{\phi}\frac{1}{w_{\left[\tfrac{n}{2}\right]}^{s_{\phi}}} \ {\rm for }\  |n|> 2m_{\phi} , 
\end{equation}
where the odd or even terms of $\bar b$ are equal to $\bar b^{*}$ and the others are set to zero.
\end{remark}

\begin{remark}\label{rmk:notation}
For the remainder of the paper the subscripted constants $s_\phi $, $m_\phi $, $M_\phi $, and $r_\phi $ introduced in the previous remark will be kept fixed. As indicated, they refer to the parameters associated with the 
computation of the bound states. In the next section (see eq. \eqref{constDef}) a new set of parameters 
$s$, $m $, $M$ will be chosen for the computation of the eigenvalues and eigenfunctions of the
linearization. The purpose of adopting the subscript  notation for $s_\phi $, $m_\phi $ etc. is to avoid confusion of the two sets of parameters. 
   
\end{remark}
\subsection{Solving the eigenvalue problem}\label{sec:eigenvalue}
\label{CompEig}
Once the solution $\phi(x)$ is computed,  the eigenvalue problem consists in solving \eqref{c},\eqref{d} for the unknowns $\beta, w(x), z(x)$ 
 As before, we expand $w(x)$ and $z(x)$ w.r.t. the Fourier-sine basis
$$
w(x)=\sqrt 2\sum_{n\geq 1}c_{n} \sin(\pi n x),\quad z(x)=\sqrt 2\sum_{n\geq 1}d_{n}\sin(\pi n x),\quad c_{n},d_{n}\in \C ,
$$
so we obtain the infinite--dimensional algebraic system
\begin{equation}\label{complete_system}
f=(f_{1},f_{2}\dots)=0,\quad f_{n}=
\left[\begin{array}{l}
(\pi n)^{2}c_{n}+(\sigma\muz-\beta)c_{n}-\sigma\sum_{\ell \geq 1}F_{n,\ell}(2 c_{\ell }-d_{\ell })\\
\\
(\pi n)^{2}d_{n}+(\sigma\mu+\beta)d_{n}+\sigma\sum_{\ell \geq 1}F_{n,\ell}( c_{\ell }-2d_{\ell })
\end{array}\right] ,
\end{equation}
to be solved for  unknowns $(\beta,c_{1},c_{2},\dots,d_{1},d_{2},\dots)$. The matrix $F=\{F_{n,\ell}\}$ corresponds to the 
term $\phi(x)^{2}$ and it is given explicitly by 
\begin{equation}\label{eq:Fnl}
F_{n,\ell}=2\Big(\sum_{p+k=n+\ell} b_{k}b_{p}-\sum_{p+k=n-\ell}b_{k}b_{p} \Big)
\end{equation}
where $b_{n}$ are the coefficients in  \eqref{eq:phiinbetas}; see Appendix B.
Since the system is invariant under rescaling of eigenfunctions, we need to choose a normalization to 
obtain isolated solutions. Rather than introducing a new equation, we remove one of the unknowns. Assume that we have computed a numerical solution $\hat x=(\bar \beta, \{\bar c_{k},\bar d_{k}\}_{k=1}^{m})$ of the system \eqref{complete_system} (for $n=1,\ldots ,m$) and that $\bar c_{j_{*}}$ is the largest value of the $\bar c_{k}$'s. Then we fix the value of $c_{j_{*}}=\bar c_{j_{*}}$ and we   remove $c_{j_{*}}$ from the unknowns.
The reduced  vector of unknowns is 
$$
x= (\beta,d_{1}, c_{1},d_{1},\dots, c_{j_{*}-1},d_{j_{*}-1},d_{j_{*}},c_{j_{*}+1},d_{j_{*}+1},\dots)
$$
and, grouping $c_{k}, d_{k}$, we  write
$$
x=(x_{0},x_{1},x_{2},\dots), \quad x_{0}=\beta, x_{j_{*}}=(d_{j_{*}}),x_{k}=(c_{k},d_{k}), k\neq j_{*}.
$$
Now choose a decay rate $s$, finite-dimensional parameter $m$, and computational parameter $M$  so that 
\begin{equation} \label{constDef}
s<s_{\phi}, \quad m =  3m_{\phi},\quad  M>m+4m_{\phi}. 
\end{equation} 
Then the $s$-norm of $x$ and the corresponding Banach space are defined by 
$$
\| x\|_{s}=\sup\{|\beta|, |d_{j_{*}}|j_{*}^{s},\sup_{k\geq 1,k\neq j_{*}}\{|c_{k}|k^{s},|d_{k}|k^{s}\}\}
$$
and 
$$
X^{s}\bydef \{x: \| x\|_{s}<\infty \},
$$
respectively.  Keeping in mind that $c_{j_{*}}$ is fixed, we look for a zero of 
$$
f(x)=(f_{1},f_{2},\dots):X^{s}\to W=X^{s-2}
$$
for $s\geq 2$ and $f_{n}$ as in \eqref{complete_system}. 
Let
$$
\bar x=(\bar \beta,\bar d_{j_{*}},\{\bar c_{k},\bar d_{k}\}_{k\geq q,\neq j_{*}}^{m})
$$
be an approximate zero for $f^{(m)}(x)$ (which can be obtained by simply removing $\bar c_{j_{*}}$ from $\hat x$).
The Jacobian $Df^{(m)}=\frac{\partial f^{(m)}}{\partial x^{(m)}}(\bar x)$ is given by 
\begin{equation}
Df^{(m)}=\left[ \begin{array}{c|cc|c|c|c|cc}
&&&&&&&\\
\frac{\partial}{\partial \beta}&\frac{\partial}{\partial c_{1}}&\frac{\partial}{\partial d_{1}}&\dots&\frac{\partial}{\partial d_{j_{*}}}&\dots&\frac{\partial}{\partial c_{m}}&\frac{\partial}{\partial d_{m}}\\
&&&&&&&\\
\end{array}\right] ,
\end{equation}
where
\begin{equation}
\frac{\partial}{\partial \beta}=\left[\begin{array}{c}
-\bar c_{1}\\
\bar d_{1}\\
-\bar c_{2}\\
\bar d_{2}\\
\vdots\\
-\bar c_{m}\\
\bar d_{m}
\end{array}\right]\qquad  \frac{\partial}{\partial c_{j}},\frac{\partial}{\partial d_{j}}=\left[\begin{array}{cc}
-2\sigma F_{1,j} & \sigma F_{1,j}\\
\sigma F_{1,j}&-2\sigma F_{1,j}\\
\vdots&\vdots\\
\pi^{2}j^{2}+\sigma\mu-\bar\beta-2\sigma F_{j,j}&\sigma F_{j,j}\\
\sigma F_{j,j}&\pi^{2}j^{2}+\sigma\mu+\bar\beta-2\sigma F_{j,j}\\
\vdots&\vdots\\
-2\sigma F_{m,j} & \sigma F_{m,j}\\
\sigma F_{m,j}&-2\sigma F_{m,j}
\end{array}\right]
\end{equation}

and,  for $k>m$, 
$$
\Lambda_{k}=\frac{\partial F_{k}}{\partial (c_{k},d_{k})}(\bar x)=\left[ \begin{array}{cc}
\pi^{2}k^{2}+\sigma\mu-\bar \beta-2\sigma F_{k,k}&\sigma F_{k,k}\\
\sigma F_{k,k}&\pi^{2}k^{2}+\sigma\mu+\bar \beta-2\sigma F_{k,k}
\end{array}\right] .
$$
Hence, according to  \eqref{eq:A},\eqref{eq:T}, the operator $A$ is defined  by the infinite-dimensional matrix 
$$
A=\left[ \begin{array}{ccc|cccc}
&&&&&&\\
&A^{(m)}&&&&&\\
&&&&&&\\
\hline
&&&(\Lambda_{m+1})^{-1}&&&\\
&&&&\ddots&&\\
&&&&&(\Lambda_{k})^{-1}&\\
&&&&&&\ddots
\end{array}\right] ,\qquad A^{(m)}Df^{(m)}\approx I,
$$
and the operator $T:X\to X$ is given by 
$$
T(x)=x-Af(x).
$$
The operator $T$ is well-defined, since the operator $A$ maps $X^{s-2}$ to $ X^{s}$ and is invertible. These properties follow from the behaviour of $ \Lambda_{k}^{-1}$ for $k>m$; in Appendix B we prove that, for sufficiently large $m$, there exists 
a constant $\mathcal C_{\Lambda}(m)$ such that 
\begin{equation}\label{eq:mClambda}
\|\Lambda_{k}^{-1}\|_{\infty}\leq\frac{\mathcal C_{\Lambda}(m)}{k^{2}}\quad (\forall k>m).
\end{equation}
Hence, fixed points of $T$ correspond to zeros of $f(x)$.

\subsubsection{Construction of the bounds Y, Z}
\label{YZ}
In deriving rigorous bounds, the most difficult terms are not actually the 
nonlinear ones (given by the product $\beta c_{k}$ and $\beta d_{k}$),  but the  linear terms, such as  $\sum_{\ell \geq 1}F_{n,\ell}(2c_{\ell }-d_{\ell })$. This is because each $F_{n,\ell}$ is defined as a convolution of $b_{k}$'s; the latter, however, are the result of the prior computation of  $\phi(x)$ and, as such, are only known to lie in certain intervals. 
Therefore, in order to design a successful scheme, we need to  find sharp estimates for the terms $F_{n,\ell}$. 

Using the notation of  Remarks \ref{rmk:enclosure} and  \ref{rmk:notation}, define 
$$
\mathcal E(q)=4^{s_{\phi}}r_{\phi}^{2}\frac{\alpha^{(2)}_{q}}{w_{q}^{s_{\phi}}}+2r_{\phi}2^{s_{\phi}}\sum_{j=-2m_{\phi}}^{2m_{\phi}}|\bar b_{j}|w_{q-j}^{-s_{\phi}}
$$
$$
\tilde{\mathcal E}(q)=w_{q}^{s_{\phi}}\mathcal E(q),
$$
where $\alpha_{q}^{(2)}$ is  defined in eq. \eqref{eq:alhak2} of the appendix (note that the constant $M$ in \eqref{eq:alhak2}
is to be interpreted as $M_\phi $).

\begin{lem}\label{lem:epsq}
For any $q$ 
\begin{equation}\label{eq:boundbpbk}
\left|\sum_{p+k=q}b_{k}b_{p}\right|\leq \left|\sum_{p+k=q}\bar b_{p}\bar b_{k}\right|+ \mathcal E(q)
\end{equation}
In particular, for $|q|\geq 4m_{\phi}$ 
$$
\left|\sum_{p+k=q}b_{k}b_{p}\right|\leq\mathcal E(q)\ .
$$
\end{lem}
{\it Proof.} See Appendix B. 
\eproof
\begin{remark}
\label{rmk:E}$\mathcal E(-q)=\mathcal E(q)$ and $\tilde {\mathcal E}(-q)=\tilde{\mathcal E}(q)$. The functions $\mathcal E (q)$ and $\tilde{\mathcal E}(q)$ are decreasing in $q$, for $q\ge m = 3m_\phi $. 
\end{remark}

Defining 
$$
b_{\max}\bydef \max_{n}|\bar b_{n}|
$$
$$
|\overline F_{n,\ell}|\bydef 2\left(\left|\sum_{p_{1}+p_{2}=n+\ell}\bar b_{ p_{1}}\bar b_{ p_{2}} \right|+\left|\sum_{p_{1}+p_{2}=n-\ell}\bar b_{p_{1}}\bar b_{ p_{2}} \right|\right)
$$
we list some properties of $F=\{F_{n,\ell}\}$. 

\begin{lem}\label{lem:F}

\par 
1.\  $F$ is symmetric.

2.\ $F_{n,\ell}=0$ for $n+\ell =0$. 

3.\ For any $n,\ell\geq 1$
$$
 F_{n,\ell}\in 2\left[ \sum_{|k|\leq 2m_{\phi}}b_{k}(b_{n+\ell-k}-b_{n-\ell-k})\right]\pm 8r_{\phi}(b_{\max}+r_{\phi})\frac{2^{s_{\phi}}}{(s_{\phi}-1)m_{\phi}^{s_{\phi}-1}}.
$$

4.\ For any $n,\ell\geq 1$
\begin{equation}\label{eq:boundFkl}
 |F_{n,\ell}|\leq |\overline F_{k,\ell}|  +2\mathcal E(n+\ell)+ 2\mathcal E(n-\ell).
\end{equation}

5.\ For any $n,\ell\geq 1$ such that $|n-\ell|>4m_{\phi}$

\begin{equation}\label{eq:boundFkl2}
|F_{n,\ell}|\leq \frac{2}{(n+\ell)^{s_{\phi}}}\tilde{\mathcal E}(n+\ell)+\frac{2}{(n-\ell)^{s_{\phi}}}\tilde{\mathcal E}(n-\ell).
\end{equation}

\end{lem}

{\it Proof.}
 $1.$ Follows directly from \eqref{eq:Fnl}. 
 
 $2.$ Immediate consequence of the fact that the even or the odd elements of $\{ b_{n}\}$ are zero.

$3.$ 
$$
F_{n,\ell}=2\left[ \sum_{|k|\leq 2m_{\phi}}b_{k}(b_{n+\ell-k}-b_{n-\ell-k})\right]+2\left[ \sum_{|k|>2m_{\phi}}b_{k}(b_{n+\ell-k}-b_{n-\ell-k})\right]
$$
and so
$$
\left|F_{n,\ell}- 2\sum_{|k|\leq 2m_{\phi}}b_{k}(b_{n+\ell-k}-b_{n-\ell-k})\right|\leq 2 \sum_{|k|>2m_{\phi}}|b_{k}|(|b_{n+\ell-k}|+|b_{n-\ell-k}|).
$$
From \eqref{eq:enclosurebn} it follows $|b_{n}|\leq b_{\max}+r_{\phi}$ ($ \forall n$), hence the right hand side of the previous inequality can be bounded by $2(2(b_{\max}+r_{\phi})2^{s^{\phi}}r_{\phi}\sum_{|k|>m_{\phi}}1/|k|^{s_{\phi}})\leq 8(b_{\max}+r_{\phi})r_{\phi}\frac{2^{s_{\phi}}}{(s_{\phi}-1)m_{\phi}^{s_{\phi}-1}}$.

$4.$ Combine \eqref{eq:Fnl} and \eqref{eq:boundbpbk}.

$5.$ Since $\bar b_{n}=0$ for $|n|\geq 2m_{\phi}$, the estimates follows from \eqref{eq:boundFkl} and the definition of $\tilde{\mathcal E}(q)$.
\eproof

We define a constant $C_{\Lambda}$  satisfying  \eqref{eq:CLambda} by 
$$
 C_{\Lambda}=\frac{\mathcal C_{\Lambda}(M)}{M^{2}} , 
 $$
where $\mathcal C_{\Lambda}(M)$ has been introduced in  \eqref{eq:mClambda}.

The bounds $Y_{k}$ for $k=1,\dots,m$ are defined as in \eqref{eq:Y}. 
The next lemma provides a uniform bound for the tail  part of $f(\bar x)$.
\begin{lem}\label{lem:H}
Let $\mathcal H = \mathcal H(M)$ be  the vector
$$
\mathcal H \bydef 2\sum_{1\leq\ell\leq m} \left(  \tilde{\mathcal E}(M+\ell)+\frac{1}{(1-\tfrac{\ell}{M})^{s_{\phi}}}\tilde{\mathcal E}(M-\ell)  \right)   \left[\begin{array}{r}
|(2\bar c_{\ell }-\bar d_{\ell })|\\
|\bar c_{\ell }-2\bar d_{\ell })|
\end{array} \right].
$$ 
Then  $|f(\bar x)|\leq\frac{1}{k^{s_{\phi}}}\mathcal H(M)$ for all $k\geq M$.
\end{lem}

{\it Proof.}  See Appendix B. \eproof

 Since $s<s_{\phi}$, the previous lemma implies that $|f_{k}(\bar x)|\leq \frac{1}{w_{k}^{s}}\frac{\mathcal H(M)}{M^{s_{\phi}-s}}$. We therefore set
$$
Y_{M}\bydef C_{\Lambda}\frac{|\mathcal H(M)|_{\infty}}{M^{s_{\phi}-s}}.
$$

As for the definition of the bounds $Z_{k}$ and $Z_{M}$, it is convenient to have a formula for the coefficients $c_{k,i}$ defined in  \eqref{eq:expansion}. To this end, we write  
$Df(\bar x+ru)- J^{\dag}=$
\begin{equation*}
\left[ \begin{array}{rccccccccrllcc}
-ru_{1}&-ru_{0}&0&&0&\dots&&&|&&&&\\
ru_{1}&0&ru_{0}&&0&\dots&&&|&&\\
\vdots&&&\ddots\\
-ru_{j_{*}}&0&0&&0&\dots&&&|&&\\
ru_{j_{*}}&&&&ru_{0}&&&&|&&\\
\vdots&&&&0&\ddots&&&|&&*&*\\
-ru_{m}&&&&&&-ru_{0}&&|&&*&*\\
ru_{m}&&&&&&&ru_{0}&|&&\\
-&-&-&-&-&-&-&-&-&&\\
-ru_{m+1}&-2\sigma F_{m+1,1}&\sigma F_{m+1,1}&\dots&\sigma F_{m+1,j_{*}}&&&&&-ru_{0}&0\\
ru_{m+1}&\sigma F_{m+1,1}&-2\sigma F_{m+1,1}&\dots&-2\sigma F_{m+1,j_{*}}&&&&&0&ru_{0}\\
\vdots&\vdots&\vdots&\vdots&&&&&&&&\ddots\\
-ru_{k}&-2\sigma F_{k,1}&\sigma F_{k,1}&\dots&\sigma F_{k,j_{*}}&*&*\\
ru_{k}&\sigma F_{k,1}&-2\sigma F_{k,1}&\dots&-2\sigma F_{k,j_{*}}&*&*\\
\vdots&\vdots
\end{array}\right]
\end{equation*}
where the $(k,j)$-block is
$$
\begin{array}{cc}
*&*\\
*&*
\end{array}=\sigma \left[\begin{array}{cc}
-2F_{k,j}& F_{k,j}\\
F_{k,j}&-2F_{k,j}
\end{array}\right].
$$
Thus we have
$$
c_{\ \cdot,1}=\sigma \left[
\begin{array}{c}
{\displaystyle -2\sum_{j=m+1}^{\infty}F_{1,j}v_{j}+\sum_{j=m+1}^{\infty}F_{1,j}v_{j}}\\
{\displaystyle \sum_{j=m+1}^{\infty}F_{1,j}v_{j}-2\sum_{j=m+1}^{\infty}F_{1,j}v_{j}}\\
\vdots\\
{\displaystyle -2\sum_{j=m+1}^{\infty}F_{m,j}v_{j}+\sum_{j=m+1}^{\infty}F_{m,j}v_{j}}\\
{\displaystyle \sum_{j=m+1}^{\infty}F_{m,j}v_{j}-2\sum_{j=m+1}^{\infty}F_{m,j}v_{j}}\\
{\displaystyle F_{m+1,j_{*}}v_{j_{*}}-2\sum_{j=1,j\neq j_{*}, m+1}^{\infty}F_{m+1,j}v_{j}+\sum_{j=1,j\neq j_{*},m+1}^{\infty}F_{m+1,j}v_{j}}\\
{\displaystyle -2F_{m+1,j_{*}}v_{j_{*}}+\sum_{j=1,j\neq j_{*}, m+1}^{\infty}F_{m+1,j}v_{j}-2\sum_{j=1,j\neq j_{*},m+1}^{\infty}F_{m+1,j}v_{j}}\\
\vdots\\
{\displaystyle F_{k,j_{*}}v_{2}-2\sum_{j=1,j\neq j_{*},k}^{\infty}F_{k,j}v_{j}+\sum_{j=1,j\neq j_{*},k}^{\infty}F_{k,j}v_{j}}\\
{\displaystyle -2F_{k,j_{*}}v_{2}+\sum_{j=1,j\neq j_{*},k}^{\infty}F_{k,j}v_{j}-2\sum_{j=1,j\neq j_{*},k}^{\infty}F_{k,j}v_{j}}\\
\vdots
\end{array}
\right],\quad c_{\ \cdot,2}=\left[ \begin{array}{c}
-u_{1}v_{0}-u_{0}v_{1}\\
 u_{1}v_{0}+u_{0}v_{1}\\
 \vdots\\
 -u_{j_{*}}v_{0}\\
u_{j_{*}}v_{0}+u_{0}v_{j_{*}}\\
\vdots\\
-u_{k}v_{0}-u_{0}v_{k}\\
u_{k}v_{0}+u_{0}v_{k}\\
\vdots
\end{array}\right]
$$
Since $|u_{j}|,|v_{j}|\leq j^{-s}$, we obtain the estimates
$$
|c_{k,1}|\leq\left[
\begin{array}{c}
3{\displaystyle \sum_{j=m+1}^{\infty}|F_{k,j}|j^{-s}}\\
3{\displaystyle \sum_{j=m+1}^{\infty}|F_{k,j}|j^{-s}}
\end{array}
\right]\quad k=1,\dots,m,\quad |c_{k,1}|\leq\left[
\begin{array}{c}
3{\displaystyle \sum_{j=1,j\neq k}^{\infty}|F_{k,j}|j^{-s}}-2|F_{k,j_{*}}|j_{*}^{-s}\\
3{\displaystyle \sum_{j=1,j\neq k}^{\infty}|F_{k,j}|j^{-s}}-|F_{k,j_{*}}|j_{*}^{-s}
\end{array}
\right]\quad k\geq m+1
$$
and
$$
|c_{j_{*},2}|\leq j_{*}^{-s}\left[\begin{array}{c}
2\\
4
\end{array}\right]\quad|c_{k,2}|\leq 4k^{-s}\left[\begin{array}{c}
1\\
1
\end{array}\right]\quad k\geq1, k\neq j_{*} ,
$$
where we used the fact that $|z_{1}z_{2}|=\max\{|Re(z_{1}z_{2}|,|Im(z_{1}z_{2}| \}\leq 2|z_{1}||z_{2}|$.

Note that the vectors  $c_{k,1}$ are given as series. We can provide a bound  using  formula \eqref{eq:boundFkl}. 

Define
\begin{equation*}
H^{1}(k)\bydef3\left\{\begin{array}{ll}
{\displaystyle \sum_{j=m+1}^{4m_{\phi}+k} |\bar F_{k,j}|j^{-s}+2\frac{\mathcal  E(m+1+k)+\mathcal E(m+1-k)}{(s-1)m^{s-1}}} & k\leq m\\
\begin{array}{l}
{\displaystyle \sum_{j=\max\{1,k-4m_{\phi}\},j\neq k}^{k+4m_{\phi}} |\bar F_{k,j}|j^{-s}+ 2\sum_{j=1}^{k-1}\mathcal E(j-k)j^{-s}}\\
{\displaystyle +\mathcal E(1)\frac{2}{(k+1)^{s}}+\mathcal E(2)\frac{2}{(s-1)(k+1)^{s-1}}+2\mathcal E(k+1)+\mathcal E(k+2)\frac{2}{(s-1)}}
 \end{array}
 &  k>m \ .
\end{array}\right.
\end{equation*}
Then we have 
\begin{lem}\label{lem:Z1}
$$
|c_{k,1}|_{\infty}\leq H^{1}(k), \quad k\geq 1.
$$
\end{lem}
{\it Proof.} See  Appendix B. \eproof 

In view of Lemma \ref{lem:Z1},  we define  $Z^{1}$, $Z^{2}$ as the vectors with components
$$
Z^{1}_{k}=H^{1}(k)\left[\begin{array}{c}
1\\
1
\end{array}\right] ,\quad k=1,\dots ,M-1,$$
$$
Z^{2}_{k}=4k^{-s}\left[\begin{array}{c}
1\\
1
\end{array}\right]\ {\rm for}\  k\neq j_{*},  \ Z^{2}_{j_{*}}=2j_{*}^{-s}\left[\begin{array}{c}
1\\
2
\end{array}\right] ,\quad k=1,\dots ,M-1.
$$

The final pieces are the bounds $Z^{1}_{M}, Z^{2}_{M}$ satisfying \eqref{eq:unifest} that will give the tail bound $Z_{M}$.  
Clearly, we can set $ Z^{2}_{M}\bydef 4$, while $Z^{1}_{M}$ has to be defined as a uniform bound (up to $w_{k}^{-s}$)  of $H^{1}(k)$ for $k\geq M$.
\begin{lem}\label{lem:Z1M}
Define
$$
Z^{1}_{M}\bydef 
6\left[\sum_{p= 1}^{4m_{\phi}} \left( \left| \sum_{p_{1}+p_{2}=p}\bar b_{p_{1}}\bar b_{p_{2}}\right|\cdot \frac{1}{(1-\tfrac{p}{M})^{s}+1}\right)+\tilde{\mathcal  E}(M+1)+\frac{\tilde{\mathcal  E}(M+2)}{s-1}+\tilde{\mathcal  E}(1)\gamma_{M}+\tilde{\mathcal  E}(1)+\frac{\tilde{\mathcal  E}(2)}{s-1}\right]
 $$ where  
 $
 \gamma_{k}
 $ is given in \eqref{eq:gammak}.
 Then $|c_{k,1}|_{\infty}\leq \frac{1}{w_{k}^{s}}Z^{1}_{M}$, for all $k\geq M$.
 \end{lem}
 {\it Proof}.\ See Appendix B. \eproof

\subsection{The eigenvalues are simple}
\label{CompSimple}
Let $\beta$ and $W={w\choose z}$ be solution of the eigenvalue problem \eqref{c},\eqref{d}. To show that $\beta$ is simple it will be verified that there is no eigenfunction $V$ of the operator $\Nop$ orthogonal  to $W$. Define the operators  $L_{\beta}$ and $G$ by 
 $L_{\beta}(V)=(\Nop-\beta I)V$ and  
$$
G(\lambda_{0},V)=\left[ \begin{array}{c}
<W,V>\\
\lambda_{0}W+L_{\beta}(V)
\end{array}\right],\quad \lambda_{0}\in\R, V\in C^{2}_{0}([0,1]),
$$
respectively. 
\begin{lem} \label{lem:Simple}
If $X_{0}=(0,0)$  is a  locally unique solution of $G(X)=0$, then the eigenvalue $\beta $ is simple.
\end{lem}
\textit{Proof.} Assume that $\lambda_{0}=0$, $V=0$ is a locally unique solution, 
but that $\beta$ is not simple. Then there exists a function $\mathcal V$ such that $<W,\mathcal V>\,=0$ and $L_{\beta}(\mathcal V)=0$. However, this implies that $X_{\lambda}=(0,\lambda\mathcal V)$ is a solution for every $\lambda\in \R$,  so
the zero--solution is not locally unique. Contradiction. \eproof

To apply the spectral method of Section \ref{sec:eigenvalue}, we recast $G(\lambda_{0},V)=0$ as an infinite--dimensional algebraic system with unknowns $x=(\lambda_{0},\{ c_{n}\},\{d_{n}\})$. Suppose that $(\beta, \{ c'_{n},d'_{n}\}_{n\geq 1})$ represents the eigenpair $(\beta,W)$, that is  $f(\beta,\{ c'_{n},d'_{n}\}_{n\geq 1})=0$, see \eqref{complete_system}.  Then we introduce the system 
$$
g=(g_{0},g_{1},g_{2},\dots)(x) =0 
$$
given by
$$
g_{0}=\sum_{n\geq 1}(c'_{n}c_{n}+d'_{n}d_{n})
\quad {\rm and} \quad 
g_{n}=\lambda_{0}\left[ \begin{array}{l}
c'_{n}\\
d'_{n}
\end{array}\right]+f_{n}(\beta,\{c_{n}\}, \{d_{n}\}),\quad n\geq 1.
$$
We adapt the radii--polynomial technique to  check that the zero--solution is locally unique. The construction of the fixed point operator and of the bounds are very similar to Section~\ref{YZ} and hence omited.

Let a numerical approximate solution $\bar x\approx 0$ be given. Then, if the computation results in a radius $r$ so that $0\in B_{\bar x}(r)$, we conclude that  $x=0$ \textit{is the}  locally unique solution and $\beta$ is simple.

\begin{remark}
i) The operator $g$ is linear in $x$, therefore the radii polynomial have degree one. 

ii) The introduction of the unknown $\lambda_{0}$ is technical; its purpose is to   {\it balance} the number of equations with the number of unknowns. However, since the operator $\Nop$ has no generalized eigenvectors, the system $g(x)=0$ cannot have any solutions with $\lambda_{0}\neq 0$. As a result,  Lemma~\ref{lem:Simple} could be rephrased to  say that $X_{0}=(0,0)$  is a  locally unique solution of $G(X)=0$ if and only if the eigenvalue $\beta $ is simple.
\end{remark}

\section{Numerical results}
\label{Num}
\subsection{Checking $\Gamma$}
Suppose the Fourier coefficients $c_{n}$, $d_{n}$ of $w(x)$ and $z(x)$ have been proved to be in a ball of radius $r$ in the space $X^{s}$ around the numerical approximation $\bar c_{n}, \bar d_{n}$. This means that 
$$
|c_{n}-\bar {c}_{n}|\leq \frac{r}{w_{n}^{s}},\qquad |d_{n}-\bar {d}_{n}|\leq \frac{r}{w_{n}^{s}}, \quad \forall n\geq 1. 
$$
We can then finally check condition (A) of the introduction; i.e. we verify that the $\Gamma$--coefficients are bounded away from zero. It can be shown \cite{BLT12} that  
$\Gamma \propto [{\Psi ^{(2)}}]'(1) \propto \sum_{n\geq 1}(-1)^{n}n(c_{n}-d_{n}) $
(where ``$\propto $" means ``proportional"), so  the enclosure of the Fourier coefficients implies
\begin{equation}\label{eq:gamma-proof}
\begin{split}
\Gamma \propto\sum_{n\geq 1}(-1)^{n}n(c_{n}-d_{n})&\in\sum_{n=1}^{m}(-1)^{n}n(\bar c_{n}-\bar d_{n})\pm r\sum_{n=1}^{\infty}\frac{1}{n^{s-1}}(1+i)\\
&\in\sum_{n=1}^{m}(-1)^{n}n(\bar c_{n}-\bar d_{n})\pm r\left(1+\frac{1}{s-2}\right)(1+i) ,
\end{split}
\end{equation}
Thus, if zero does not belong to the set  on the right hand side of \eqref{eq:gamma-proof},  $\Gamma $ does not vanish.
\subsection{Results}
We now describe some of the computational results obtained by the method discussed above. The results are rigorous, 
since all computations are performed in interval arithmetics.
\subsubsection{Bounded state solution of the NLS}
Table~\ref{tab1} below shows the results for the ground state, as well as the first and second exited states --  i.e. solutions $\phi(x)$ of  \eqref{groundA1},\eqref{groundB1} with $j=0,1,2$ --
for the focusing ($\sigma=1$, left half) and defocusing ($\sigma=-1$, right half) cases
and three  values of the chemical potential $\mu$. The numerical solution $\bar b$ of the Galerkin projection    $f^{m}(b)=0$ was computed by the Newton method to accuracy $|f^{m}(\bar b)|<10^{-13}$. The table lists the finite-dimensional parameter $m_{\phi}$, the decay-rate parameter $s_{\phi}$ and the radius of the ball in the space $X^{s_{\phi}}$ around the numerical solution $\bar b$ within which the solution of the infinite-dimensional problem is guaranteed to exist.
\begin{table}
$$
\begin{array}{cccccl}
\sigma  & Nodes  & \mu & m_{\phi}& s_{\phi} &\qquad r_{\phi}   \\
\hline
1  & 0  & 12.898 & 18 &4 & 4.0089 \cdot 10^{-13}  \\
1  & 0  & 43.273 & 24 & 4 & 4.7045 \cdot 10^{-10}  \\
1& 0& 80.518& 30&3.5&2.9894 \cdot 10^{-9}\\
1 & 1 & 12.898& 30& 4& 7.6398 \cdot 10^{-13}\\
1 & 1& 43.273 & 30 & 3.5 & 1.3889 \cdot 10^{-11}\\
1 & 1& 80.518 & 30 &3.2 & 2.5216 \cdot 10^{-8}\\
1& 2 & 12.898 & 44 & 3.1 & 5.5127 \cdot 10^{-12}\\
1 & 2 & 43.273 & 58 &3.1 & 1.0442\cdot 10^{-11}\\
1&2&80.518 & 80 & 3&  3.7820 \cdot 10^{-14}
\end{array}
\quad
\begin{array}{cccccl}
\sigma  & Nodes  & \mu & m_{\phi}& s_{\phi} &\qquad r_{\phi}   \\
\hline
-1  & 0  & 89.237 & 20 & 4 & 9.0256 \cdot 10^{-8}  \\
-1& 0& 161.521& 30&3.5 & 3.1080 \cdot 10^{-9}\\
-1 & 0 & 254.916 & 36& 3.1 & 6.7484\cdot 10^{-9}\\
-1& 1 & 89.237 & 20 & 4 & 4.9679\cdot 10^{-8}\\
-1& 1& 161.521& 30&3.1 & 8.4114 \cdot 10^{-10}\\
-1 & 1 & 254.916 & 34& 3.1 & 4.0952\cdot 10^{-8}\\
-1& 2 & 89.237 & 16 & 4& 2.2678\cdot 10^{-12}\\
-1&2& 161.521 & 54 & 3.5 & 1.2129 \cdot 10^{-11}\\
-1&2& 254.916 & 54 & 3& 1.0558 \cdot 10^{-14}\\
\end{array}
$$
\caption{Enclosure of bound states}
\label{tab1}
\end{table}
\begin{figure}[htbp]
\begin{center}
\includegraphics[width=1 \textwidth, height=270 pt]{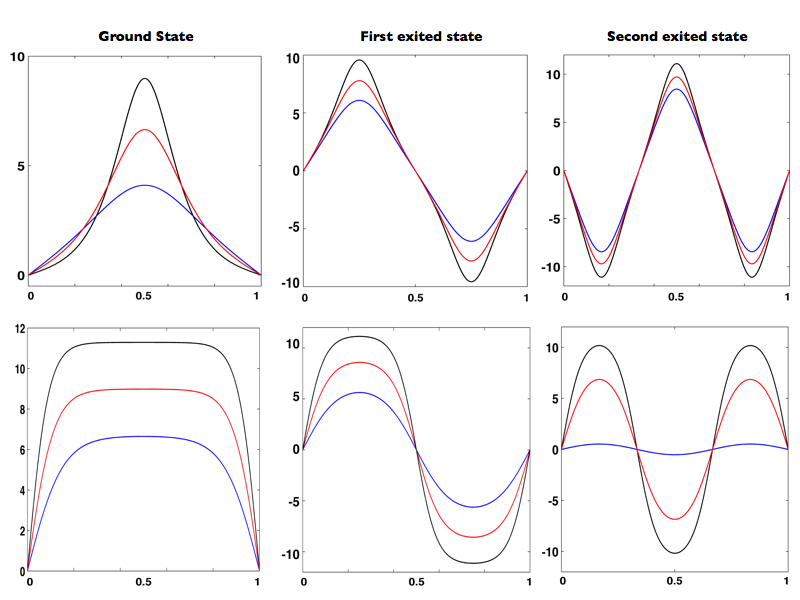}
\caption{Bounded states for three values of $\mu $ (increasing blue $\to $ red $\to $ black). Top row: focusing case
($\sigma =1$); bottom row: defocusing case ($\sigma =-1$). }
\label{default}
\end{center}
\end{figure}
\subsubsection{Enclosure of spectra and $\Gamma$--values}
For a given value of $\mu$, we considered the three different bounded states ($j=0,1,2$) computed previously.  
Representative data (two values of $\mu $; one focusing, one defocusing) for the first three (non-zero) eigenvalues and the corresponding $\Gamma $--values are listed in Tables~\ref{tabu2} and \ref{tabu3}: $r$ denotes the radius of the ball in the space $X^{s}$ (around the approximate eigenvalue and associated eigenfunction)  within which the ``true" solution $(\beta,\{c_{k},d_{k} \})$ of \eqref{complete_system} is guaranteed to exist.  The last column contains the enclosure intervals  of the corresponding $\Gamma$-values.
\begin{table}
$$
\begin{array}{c|c|c|c|c}
& Eigenvalue & s & r & \Gamma \in\\
\hline
& 13.413 & 3&  1.189\cdot 10^{-7}   &  0.9771 \pm 2.379\cdot 10^{-7}    \\
  Ground\  State &238.868  & 3& 0.355 \cdot 10^{-7}     & -4.6741  \pm   0.709  \cdot 10^{-7}   \\
   &791.201 & 3&     0.622\cdot 10^{-7} &     -8.8452 \pm  1.243\cdot 10^{-7}\\
   \hline
 &  90.461 &  3&   2.027\cdot 10^{-8} &   3.4138 \pm   4.0546\cdot 10^{-8}\\
1^{st } \ Exited\  State  & 426.79 &   3&  1.842\cdot 10^{-8} & -6.5821\pm   3.685\cdot 10^{-8}\\
  &  743.05&   3&  2.264\cdot 10^{-8} &-8.6971\pm  4.5283\cdot 10^{-8} \\
  &40.30\pm 15.51i& 3& 6.258\cdot 10^{-8}  & (-0.4929\pm 1.3720i) \pm (1+i)\cdot 1.2517 \cdot 10^{-7} \\
  \hline
  &221.73 &3& 1.357\cdot 10^{-9} & 5.462\pm 4.763\cdot 10^{-9}\\
   &676.54 &3& 1.821\cdot 10^{-9} & -8.4667 \pm 5.003\cdot 10^{-9}   \\
 2^{nd } Exited \ State &59.95\pm 25.55 i &3& 2.932 \cdot 10^{-9}& (0.6855\mp 1.5570i) \pm (1+i)\cdot 5.863\cdot 10^{-9}\\
  &120.36\pm 33.13i&3& 2.402\cdot 10^{-9} &(0.4625\mp 2.8174 i)\pm (1+i)\cdot4.804\cdot 10^{-9}  \\
\end{array} 
$$
\vspace{-3ex}
\caption{Eigenvalues and $\Gamma ${\it s}: {\bf focusing case}, $\mu=43.273$ }
\label{tabu2}
\end{table}

\begin{table}
$$
\begin{array}{c|c|c|c|c}
& Eigenvalue & s & r & \Gamma \in\\
\hline
&78.671 &3&   5.1339\cdot 10^{-6}&   1.7575\pm   1.0268\cdot 10^{-5}\\
Ground \  State&360.29 &3&    2.0547\cdot 10^{-6}&  -5.7589\pm  4.1094\cdot 10^{-6}\\
&943.45 &  3& 3.3776\cdot 10^{-6} & -9.8213\pm  6.7551\cdot 10^{-6}\\
\hline
&5.1026&   3&   2.2796\cdot 10^{-4} & 3.7268\cdot 10^{-3}\pm    4.5592\cdot 10^{-4} \\
1^{st}\ Exited \ State &284.60 &  3&   1.0601\cdot 10^{-5} &  -5.1979\pm   2.1203\cdot 10^{-5} \\
&861.30&  3&9.8321\cdot 10^{-6}&   -9.6062\pm   1.9664\cdot 10^{-5} \\
\hline
&24.184 & 2.8& 5.277\cdot 10^{-12}& -0.130\pm 2.356\cdot 10^{-11}\\
2^{nd }\ Exited\  State&452.93& 2.8& 1.176\cdot 10^{-12}&-7.229\pm 3.993\cdot 10^{-12} \\
&774.05 & 2.8& 1.369\cdot 10^{-12}&-9.397\pm 4.067\cdot 10^{-12}\\  
\end{array} 
$$
\vspace{-3ex}
\caption{Eigenvalues and $\Gamma ${\it s}: {\bf defocusing case}, $\mu=254.916$}
\label{tabu3}
\end{table}

\begin{figure}[htbp]
\begin{center}
\includegraphics[scale=0.7]{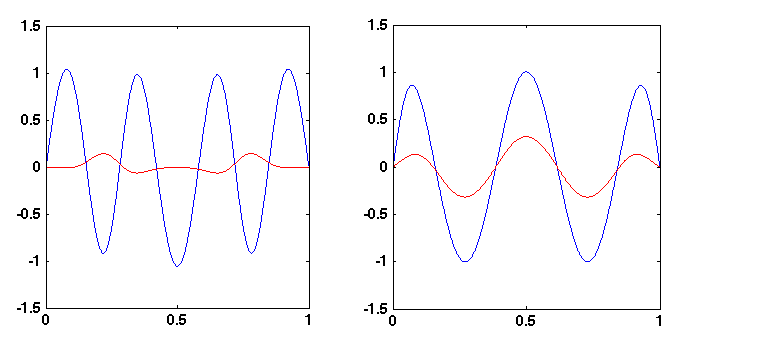}
\caption{Eigenfunctions $w(x)$ (blue) and $z(x)$ (red); cf. equations \eqref{c},\eqref{d}. Left panel: eigenfunction associated with eigenvalue $\beta=743.05$ for the first exited state with $\sigma=1$ (focusing). Right panel: $\beta=360.29$ for the ground state and $\sigma =-1$ (defocusing).
 }
\label{default}
\end{center}
\end{figure}

\begin{remark}
It may be surprising and/or confusing that the eigenvalues and $\Gamma $ values as well as their enclosure intervals are sometimes written in real form and sometimes written in complex form. 
To explain this,  we first note that all computations are carried out in complex Banach spaces. However, if the numerical (approximate) solution (i.e the centre of the enclosure interval) is real, then  
the \textit{exact} solution is real as well. Indeed,  if the (exact) solution was complex, its complex conjugate would be a solution as well, which would fall in the same 
enclosure ball. This is impossible by  uniqueness. 
\end{remark}

\section{Concluding remarks}
\label{Dis}
In this paper we analyzed important aspects of a realistic
model for a one-dimensional BEC by numerical means. Since the 
results are derived from a computational scheme that is based on the radii--polynomial technique in conjunction with interval arithmetic, they are mathematically rigorous and can be used to complement and complete analytical proofs, such as the
controllability proof given in \cite{BLT12}.  The method adopted is general and flexible; as a result,  both the focusing and defocusing cases as well as ground and excited states
can be treated within the same computational framework. 
Specifically, we 
\begin{enumerate}
\item rigorously computed the ground \textrm{and} (the first two\footnote{This number is completely arbitrary; there is no restriction in principle to rigorously 
determining any number of bound states.  A similar comment applies to the number of eigenvalues of the linearization.}) excited states;
\item rigorously computed finitely many eigenvalues and eigenfunctions of the linearization (around the bound states determined previously); 
\item proved (by rigorous numerics) that the eigenvalues are simple (B); 
\item rigorously verified the ``$\Gamma$--condition" (A).  
\end{enumerate}
The model studied in this paper has considerable interest in its stated form, both from the physical and the  mathematical
point of view (as for the latter, we note that only very few applications of  rigorous numerics to \textit{infinite}-dimensional problems 
exist to date). However, there are some obvious generalizations that immediately come to mind, such as
the whole-space problem (with a suitable potential, such as the harmonic oscillator) to replace Dirichlet boundary conditions 
and/or higher space dimensions. These generalizations 
are subject to current research by the authors and will be reported on 
 in the future.

Furthermore, in addition to presenting a study of an (important) particular model, we also view this paper as a case study 
that illustrates the general utility and flexibility of the rigorous--numerics paradigm.   
We believe that the latter will find applications with other important problems in mathematics and science and will thus  become 
a valuable tool in the arsenal of mathematicians,  physicists, and scientists at large.  
\section*{Acknowledgements} This work was partially supported by grant {MTM2011-24766} of the MICINN, Spain, and the \textit{Natural Sciences and 
Engineering Research Council of Canada (NSERC)}.   
This research was begun during a three-month stay of the second author at the \textit{Basque Center for Applied Mathematics (BCAM)},
who would like to thank  BCAM for its hospitality and financial support.

\section{Appendix}

The integral
\begin{equation}
\begin{array}{c}\label{eq:intsenis}
\int_{0}^{1}\sin(\pi k x)\sin(\pi l x)\sin(\pi p x)\sin(\pi n x)dx=\\
\frac{1}{8}\Big(\delta_{p+k-l-n}-\delta_{p+k-l+n}+\delta_{p-k+l-n}-\delta_{p-k+l+n}-\delta_{p+k+l-n}+\delta_{p+k+l+n}-\delta_{p-k-l-n}+\delta_{p-k-l+n}\Big)
\end{array}
\end{equation}
is readily computed. Hence, given $\phi(x)=\sum_{n\in \Z}b_{n}\sin(n\pi x)$, we have
\begin{equation}
\begin{split}
<\phi^{3},\sqrt 2\sin(\pi n \bullet)>&=4\int_{0}^{1}\Big[\sum_{k\in \Z}b_{k}\sin(\pi k x)\sum_{\ell \in \Z}b_{\ell }\sin(\pi l x)\sum_{p\in \Z}b_{p}\sin(\pi p x)\Big]\sin(\pi n x)dx\\
&=4\sum_{k,\ell,p\in \Z}b_{k}b_{\ell }b_{p}\int_{0}^{1}\sin(\pi k x)\sin(\pi l x)\sin(\pi p x)\sin(\pi n x)dx\ .\\
\end{split}
\end{equation}
Using \eqref{eq:intsenis} and the property $b_{-k}=-b_{k}$,  eq. \eqref{eq:cubica} follows. 
\subsection{Appendix A: analytical estimates for the enclosure of $\phi(x)$ }\label{AppendixA}
\subsubsection*{Bounds}
The definition of the bounds $Y,Z$ is the same as in \cite{MR2718657}, so we refer to that paper for a detailed explanation. 
We first recall the definition of some constants:
\begin{equation}\label{eq:gammak}
\gamma_{k}=2\left[\frac{k}{k-1} \right]^{s}+\left[ \frac{4\ln(k-2)}{k}+\frac{\pi^{2}-6}{3}\right]\left[\frac{2}{k} +\frac{1}{2}\right]^{s-2}
\end{equation}
\begin{equation}\label{eq:alhak2}
\alpha_{k}^{(2)}=\left\{\begin{array}{ll}
4+\frac{1}{2^{2s-1}(2s-1)} & k=0\\
2[2+\tfrac{1}{2^{s}}+\tfrac{1}{3^{s}}+\frac{1}{3^{s-1}(s-1)}]+\sum_{k_{1}=1}^{k-1}\frac{k^{s}}{k_{1}^{s}(k-k1)^{s}} &1\leq k\leq M-1\\
2[2+\tfrac{1}{2^{s}}+\tfrac{1}{3^{s}}+\frac{1}{3^{s-1}(s-1)}]+\gamma_{k}& k\geq M
\end{array}\right.
\end{equation}
$$
\alpha_{k}^{(3)}=\left\{\begin{array}{ll}
\alpha_{0}^{(2)}+2\sum_{k_{1}=1}^{M-1}\frac{\alpha_{k_{1}}^{(2)}}{k_{1}^{2s}}+\frac{2\alpha_{M}^{(2)}}{(M-1)^{2s-1}(2s-1)}& k=0\\
\sum_{k_{1}=1}^{M-k-1}\frac{\alpha_{k_{1}+k}^{(2)}k^{s}}{k_{1}^{s}(k+k_{1})^{s}}+\alpha_{M}^{(2)}k^{s}[\frac{1}{(M-k)^{s}M^{s}}+\frac{1}{(M-k)^{s-1}M^{s}(s-1)}]&\\
+\alpha_{k}^{(2)}+\sum_{k_{1}=1}^{k-1}\frac{\alpha_{k_{1}}^{(2)}k^{s}}{k_{1}^{s}(k-k_{1})^{s}}+\alpha_{0}^{(2)}+\sum_{k_{1}=1}^{M-1}\frac{\alpha_{k_{1}}^{(2)}k^{s}}{k_{1}^{s}(k+k_{1})^{s}}+\frac{\alpha_{M}^{(2)}}{(M-1)^{s-1}(s-1)} & 1\leq k\leq M-1\\
\alpha_{M}^{(2)}[2+\tfrac{1}{2^{s}}+\tfrac{1}{3^{s}}+\frac{1}{3^{s-1}(s-1)+\frac{1}{(M-1)^{s-1}(s-1)}}+\gamma_{k}]&\\
+\alpha_{0}^{(2)}+\sum_{k_{1}=1}^{M-1}(\frac{\alpha_{k_{1}}^{(2)}}{k_{1}^{s}}[1+\frac{M^{s}}{(M-k_{1})^{s}}]) & k\geq M
\end{array}\right.
$$

$$
\varepsilon^{(3)}_{k}=\frac{2\alpha_{M}^{(2)}}{(s-1)(M-1)^{s-1}(M+k)^{s}}+\sum_{k_{1}=M}^{M+k-1}\frac{\alpha^{(2)}_{k_{1}-k}}{w_{k_{1}}^{s}w_{k_{1}-k}^{s}}
$$
$$
\tilde\alpha^{(3)}_{M}\bydef\max\{\alpha_{k}^{(3)}: k=0\dots M \}.
$$
Define the bound $Y$ as
\begin{equation}
Y_{k}\bydef\left\{\begin{array}{ll}
|[A^{(m)}f^{m}(\bar b)]_{k}|,\quad &k=1,\dots,m\\
|\Lambda_{k}^{-1}f_{k}(\bar b)| &k=m+1,\dots,M-1\\
0&k\geq M
\end{array}\right. .
\end{equation}
For the bound $Z_{k}$ we first define 
\begin{equation}
 (Z^{0})_{k}=\left\{\begin{array}{ll}
\left[|I-A^{m}Df^{(m)}|\{w_{j}^{-s}\}_{j=1}^{m}\right]_{k}, \quad &k=1,\dots,m\\
~ 0, & k>m
\end{array}\right.\ .
\end{equation} 
and
\begin{equation}
 (Z^{1})_{k}=\left\{\begin{array}{ll}
{\displaystyle \sum_{\substack{
k_{1}+k_{2}+k_{3}=k \\
|k_{1}|,|k_{2}|< m,m\leq|k_{3}|<M}}|\bar b_{k_{1}}||\bar b_{k_{2}}|\tfrac{1}{w_{k_{3}}^{s}}} +\|\bar b\|_{s}^{2}\varepsilon^{(3)}_{k} \quad &k=1,\dots,m\\
\\
{\displaystyle \sum_{\substack{
k_{1}+k_{2}+k_{3}=k \\
|k_{1}|,|k_{2}|< m,|k_{3}|<M}}|\bar b_{k_{1}}||\bar b_{k_{2}}|\tfrac{1}{w_{k_{3}}^{s}}} + \|\bar b\|_{s}^{2}\varepsilon^{(3)}_{k} \quad & m+1\leq k < M
\end{array}\right.\ .
\end{equation} 
and, for any $1\leq k<M$, 
\begin{equation}
(Z^{2})_{k}\bydef\sum_{\substack{
k_{1}+k_{2}+k_{3}=k \\
|k_{1}|< m,|k_{2}|,|k_{3}|<M}}|\bar b_{k_{1}}|\tfrac{1}{w_{k_{2}}^{s}}\tfrac{1}{w_{k_{3}}^{s}} +2\|\bar b\|_{s}\varepsilon^{(3)}_{k}
\end{equation} 
\begin{equation}
(Z^{3})_{k}\bydef\sum_{\substack{
k_{1}+k_{2}+k_{3}=k \\
|k_{_{j}}|< M}}\tfrac{1}{w_{k_{3}}^{s}}\tfrac{1}{w_{k_{2}}^{s}}\tfrac{1}{w_{k_{3}}^{s}} +3\varepsilon^{(3)}_{k}\ .
\end{equation} 
Collecting all the terms, we have
\begin{equation}
 Z_{k}\bydef\left\{\begin{array}{ll}
6[|A^{m}|((Z^{1})^{m}r+2(Z^{2})^{m}r^{2}+(Z^{3})^{m}r^{3}]_{k}+Z^{0}_{k}r&1\leq k<m\\
6|\mu_{k}^{-1}|(Z^{1}_{k}r+2Z^{2}_{k}r^{2}+Z^{3}_{k}r^{3})& m\leq k<M\\
\end{array}\right.\ .
\end{equation} 
The tail bound $Y_{M}$ can be set equal to zero, while 
$$
Z_{M}=6C_{\Lambda}(\|\bar b\|^{2}_{s}\tilde\alpha^{(3)}_{M}r+2\|\bar b\|_{s}\tilde\alpha^{(3)}_{M}r^{2}+\tilde\alpha^{(3)}_{M}r^{3}).
$$
Finally, the radii polynomials are
$$
p_{k}(r)=Y_{k}+(Z^{0}_{k}+6[|A^{(m)}|(Z^{1})^{m}]_{k}-1/w_{k}^{s})r+(12[|A^{(m)}|(Z^{2})^{m}]_{k})r^{2}+(6[|A^{(m)}|(Z^{3})^{m}]_{k})r^{3},\quad 1\leq k\leq m
$$
$$
p_{k}(r)=Y_{k}+\left( \frac{6Z^1_{k}}{|\mu_{k}|}-\frac{1}{w_{k}^{s}}\right)r+\frac{12Z^2_{k}}{|\mu_{k}|}r^{2}+\frac{6Z^3_{k}}{|\mu_{k}|}r^{3} ,\quad m< k<M.
$$
$$
p_{M}(r)=6C_{\Lambda}\tilde\alpha^{(3)}_{M}r^{2}+12 C_{\Lambda}\|\bar \beta\|_{s}\tilde\alpha^{(3)}_{M}r+6C_{\Lambda}\|\bar\beta\|^{2}_{s}\tilde\alpha^{(3)}_{M}-1
$$
\subsection{Appendix B}
\subsubsection*{Construction of the matrix $F$}
The matrix $F=\{F_{n,\ell}\}$ is defined by
$$
<\phi^{2}y,\sqrt{2}\sin(\pi n\bullet)>\,=\sum_{\ell\geq 1}F_{n,\ell}\xi_{\ell}
$$ 
where $y=\sqrt{2}\sum_{k\geq 1}\xi_{k}\sin(\pi k x)$ and $\phi(x)=\sqrt 2\sum_{k\geq 1}\alpha_{k}\sin(\pi k x)$.
Now, by \eqref{eq:intsenis} and using the symmetry $b_{-k}=-b_{k}$, 
\begin{equation*}
\begin{array}{c}
{\displaystyle <\phi^{2}y,\sqrt{2}\sin(\pi n\bullet)>=}\\
{\displaystyle 4\sum_{p,k\in \Z,\ell\geq 1}\alpha_{p}\alpha_{k}\xi_{\ell}\int_{0}^{1}\sin(\pi k x)\sin(\pi \ell x)\sin(\pi p x)\sin(\pi n x)dx}=\\
{\displaystyle \frac{1}{2}\sum_{\ell\geq 1}\left[ \sum_{p+k=\ell+n}b_{p}b_{k}-\sum_{p+k=\ell-n}b_{p}b_{k}+2\sum_{p-k=n-\ell}b_{p}b_{k}-2\sum_{p-k=\ell+n}b_{p}b_{k}-\sum_{p+k=n-\ell}b_{p}b_{k}+\sum_{p+k=-n-\ell}b_{p}b_{k}\right]}=\\
{\displaystyle \frac{1}{2}\sum_{\ell\geq 1 }\xi_{\ell}\left[ 4\sum_{p+k=n+\ell}b_{p}b_{k}-4\sum_{p+k=\ell-n}b_{p}b_{k}\right]}
\end{array}
\end{equation*}
giving
$$
F_{n,\ell}=2\left[ \sum_{p+k=n+\ell}b_{p}b_{k}-\sum_{p+k=\ell-n}b_{p}b_{k}\right].
$$
\subsubsection*{Bound for $\|\Lambda^{-1}_{k}\|_{\infty} $}
Recall the definition of $\Lambda_{k}$:
$$
\Lambda_{k}=\frac{\partial F_{k}}{\partial (c_{k},d_{k})}(\bar x)=\left[ \begin{array}{cc}
\pi^{2}k^{2}+\sigma\mu-\bar \beta-2\sigma F_{k,k}&\sigma F_{k,k}\\
\sigma F_{k,k}&\pi^{2}k^{2}+\sigma\mu+\bar \beta-2\sigma F_{k,k}
\end{array}\right] .
$$
Since these are  diagonally dominated matrices, we have that
$$
\|\Lambda^{-1}\|_{\infty}\leq\max\left\{\frac{1}{|\pi^{2}k^{2}+\sigma\mu-\bar\beta-2\sigma F_{k,k}|-|F_{k,k}|}, \frac{1}{|\pi^{2}k^{2}+\sigma\mu+\bar\beta-2\sigma F_{k,k}|-|F_{k,k}|}\right\} .
$$
If $k$ is large enough, both denominators are greater than 
$$
\pi^{2}k^{2}+\sigma\mu-|\bar\beta|-3|F_{k,k}| = k^{2}\left(\pi^{2}+\frac{\sigma\mu}{k^{2}}-\frac{1}{k^{2}}(|\bar\beta|+3|F_{k,k}|)\right) .
$$
For $k>m$ and assuming $m>2m_{\phi}$,
$$
|F_{k,k}|\leq 2\left|\sum_{p_{1}+p_{2}=0}\bar b_{p_{1}}\bar b_{p_{2}}\right|+2\mathcal E(2k)+2\mathcal E(0)\leq 2\left|\sum_{p_{1}+p_{2}=0}\bar b_{p_{1}}\bar b_{p_{2}}\right|+2\mathcal E(2m)+2\mathcal E(0)=: C_{F}.
$$
Therefore for any $k>m$
$$
\|\Lambda^{-1}\|_{\infty}\leq\frac{\mathcal C_{\Lambda}(m)}{k^{2}}
$$
with
$$
\mathcal C_{\Lambda}(m)\bydef \frac{1}{\pi^{2}+\frac{\sigma\mu}{(m+1)^{2}}-\frac{1}{(m+1)^{2}}(|\bar\beta|+3C_{F})}.
$$

{\bf Proof of Lemma \ref{lem:epsq}.}

In view of \eqref{eq:enclosurebn} we have 
\begin{equation}
\begin{split}
\sum_{p+k=q}b_{p}b_{k}&\in\sum_{p+k=q}(\bar b_{p}\pm\frac{r_{\phi}2^{s_{\phi}}}{w_{p}^{s_{\phi}}})(\bar b_{k}\pm\frac{r_{\phi}2^{s_{\phi}}}{w_{k}^{s_{\phi}}})\\
&\in\sum_{p+k=q}\bar b_{p}\bar b_{k}\pm\left(2r_{\phi}2^{s_{\phi}}\sum_{p+k=q}|\bar b_{p}|w_{k}^{-s_{\phi}}+r_{\phi}^{2}4^{s_{\phi}}\sum_{p+k=q}w_{k}^{-s_{\phi}}w_{p}^{-s_{\phi}}\right)
\end{split}
\end{equation}
Then
$$
\left| \sum_{p+k=q}b_{p}b_{k}\right|\leq \left|\sum_{p+k=q}\bar b_{p}\bar b_{k}\right| +\left(2r_{\phi}2^{s_{\phi}}\sum_{p+k=q}|\bar b_{p}|w_{k}^{-s_{\phi}}+r_{\phi}^{2}4^{s_{\phi}}\sum_{p+k=q}w_{k}^{-s_{\phi}}w_{p}^{-s_{\phi}}\right).
$$
Using 
$$
\sum_{k_{1}+k_{2}=q}\frac{1}{w_{k_{1}}^{s_{\phi}}}\frac{1}{w_{k_{2}}^{s_{\phi}}}\leq\frac{\alpha_{q}^{(2)}}{w_{q}^{s_{\phi}}}
$$
\cite[Lemma A.3]{MR2718657}, gives the first assertion.
Moreover, since $\bar b_{k}=0$ for $|k|>2m_{\phi}$, the first sum is equal to zero whenever $|q|>4m_{\phi}$. \eproof 

{\bf Proof of Lemma \ref{lem:H}.}

Since $\bar x_{k}=0$ for $k>m$, we have  
\begin{equation}
f_{k}(\bar x)=\left[\begin{array}{r}
-\sum_{1\leq\ell\leq m}F_{k,\ell}(2\bar c_{\ell }-\bar d_{\ell })\\
\sum_{1\leq\ell\leq m}F_{k,\ell}(\bar c_{\ell }-2\bar d_{\ell })
\end{array} \right],\qquad \forall k>m
\end{equation}
and so
\begin{equation}
|f_{k}(\bar x)|\leq
\sum_{1\leq\ell\leq m}|F_{k,\ell}|\left[\begin{array}{r}
|(2\bar c_{\ell }-\bar d_{\ell })|\\
|\bar c_{\ell }-2\bar d_{\ell })|
\end{array} \right],\qquad \forall k>m.
\end{equation}
If $k\geq M$ we have that $k-m\geq 4m_{\phi}$; hence, by \eqref{eq:boundFkl}, 
\begin{equation}
\begin{split}
|f_{k}(\bar x)|&\leq
2\sum_{1\leq\ell\leq m}\left(\frac{1}{(k+\ell)^{s_{\phi}}}\tilde{\mathcal E}(k+\ell)+\frac{1}{(k-\ell)^{s_{\phi}}}\tilde{\mathcal E}(k-\ell) \right)  \left[\begin{array}{r}
|(2\bar c_{\ell }-\bar d_{\ell })|\\
|\bar c_{\ell }-2\bar d_{\ell })|
\end{array} \right]\\
&\leq\frac{2}{k^{s_{\phi}}}\sum_{1\leq\ell\leq m} \left(   \frac{1}{(1+\tfrac{l}{k})^{s_{\phi}}}\tilde{\mathcal E}(k+\ell)+\frac{1}{(1-\tfrac{l}{k})^{s_{\phi}}}\tilde{\mathcal E}(k-\ell)  \right)   \left[\begin{array}{r}
|(2\bar c_{\ell }-\bar d_{\ell })|\\
|\bar c_{\ell }-2\bar d_{\ell })|
\end{array} \right]\qquad \forall k\geq M.
\end{split}
\end{equation}
From the monotonicity of $\tilde{\mathcal E}(k)$ it follows that for any $k\geq M$
\begin{equation}
|f_{k}(\bar x)|\leq\frac{2}{k^{s_{\phi}}}\sum_{1\leq\ell\leq m} \left(  \tilde{\mathcal E}(M+\ell)+\frac{1}{(1-\tfrac{l}{M})^{s_{\phi}}}\tilde{\mathcal E}(M-\ell)  \right)   \left[\begin{array}{r}
|(2\bar c_{\ell }-\bar d_{\ell })|\\
|\bar c_{\ell }-2\bar d_{\ell })|
\end{array} \right]=:\frac{1}{k^{s_{\phi}}}\mathcal H(M).
\end{equation}
\eproof

{\bf Proof of Lemma \ref{lem:Z1}.}

For any $k=1,\dots,m$
\begin{equation}
\begin{split}
\sum_{j=m+1}^{\infty}|F_{k,j}|j^{-s}\leq & 2\sum_{j=m+1}^{4m_{\phi}+k}\left(\Big|\sum_{p_{1}+p_{2}=k+j}\bar p_{1}\bar p_{2} \Big|+\Big|\sum_{p_{1}+p_{2}=j-k}\bar p_{1}\bar p_{2} \Big|\right)j^{-s}\\
&+2\sum_{j=m+1}^{\infty}(\mathcal E(k+j)+\mathcal E(j-k))j^{-s}.
\end{split}
\end{equation}
Moreover, 
 \begin{equation}\label{eq:boundm1}
 \begin{split}
 \sum_{j=m+1}^{\infty}|F_{k,j}|j^{-s}&\leq\sum_{j=m+1}^{4m_{\phi}+k} |\bar F_{k,j}|j^{-s}+2\sum_{j=m+1}(\mathcal  E(m+1+k)+\mathcal E(m+1-k))j^{-s}\\
 &\leq \sum_{j=m+1}^{4m_{\phi}+k} |\bar F_{k,j}|j^{-s}+2\frac{\mathcal  E(m+1+k)+\mathcal E(m+1-k)}{(s-1)m^{s-1}}.
 \end{split}
 \end{equation}
  For $k>m$
\begin{equation}\label{eq:boundm2} 
 \begin{split}
 \sum_{j=1,j\neq k}^{\infty}|F_{k,j}|j^{-s}&\leq\sum_{j=\max\{1,k-4m_{\phi}\},j\neq k}^{k+4m_{\phi}} |\bar F_{k,j}|j^{-s}+2\sum_{j=1,j\neq k}^{\infty}(\mathcal  E(j+k)+\mathcal E(|j-k|))j^{-s}\\
 &\leq\sum_{j=\max\{1,k-4m_{\phi}\},j\neq k}^{k+4m_{\phi}} |\bar F_{k,j}|j^{-s}+2\sum_{j=1}^{k-1}\mathcal E(|j-k|)j^{-s}\\
 &\qquad \,\,+2\sum_{j=k+1}^{\infty}\mathcal E(|j-k|)j^{-s}+2\sum_{j=1}^{\infty}\mathcal  E(j+k)j^{-s}\\
 &\leq\sum_{j=\max\{1,k-4m_{\phi}\},j\neq k}^{k+4m_{\phi}} |\bar F_{k,j}|j^{-s}+ 2\sum_{j=1}^{k-1}\mathcal E(k-j|)j^{-s}\\
& +\mathcal E(1)\frac{2}{(k+1)^{s}}+\mathcal E(2)\frac{2}{(s-1)(k+1)^{s-1}}+\mathcal E(k+1)+\mathcal E(k+2)\frac{2}{(s-1)}. \quad \Box
 \end{split}
 \end{equation}

{\bf Proof of Lemma \ref{lem:Z1M}}

We need to find $Z^{1}_{M}$ such that
$$
|c_{k,1}|_{\infty}\leq \frac{1}{w_{k}^{s}}Z^{1}_M \quad k\geq M.
$$
This requires a uniform bound for $|c_{k,1}|_{\infty}$  as $k\geq M$.
First we have
\begin{equation}\label{eq:Fkj}
 \sum_{j=1,j\neq k}^{\infty}|F_{k,j}|j^{-s}\leq\sum_{j=k-4m_{\phi},j\neq k}^{k+4m_{\phi}} |\bar F_{k,j}|j^{-s}+2\sum_{j=1,j\neq k}^{\infty}(\mathcal  E(j+k)+\mathcal E(|j-k|))j^{-s}\end{equation}
 Since $M>m+4m_{\phi}$, 
 \begin{equation}\label{eq:boundZM1}
 \begin{split}
 \sum_{j=k-4m_{\phi},j\neq k}^{k+4m_{\phi}} |\bar F_{k,j}|j^{-s}&=2\sum_{p=-4m_{\phi},p\neq 0}^{4m_{\phi}} \left| \sum_{p_{1}+p_{2}=p}\bar b_{p_{1}}\bar b_{p_{2}}\right|\cdot (k-p)^{-s}\\
&=\frac{2}{k^{s}}\sum_{p= 1}^{4m_{\phi}} \left| \sum_{p_{1}+p_{2}=p}\bar b_{p_{1}}\bar b_{p_{2}}\right| \cdot \frac{1}{(1-\tfrac{p}{k})^{s}+(1+\tfrac{p}{k})^{s}}\\
&\leq \frac{2}{k^{s}}\sum_{p= 1}^{4m_{\phi}} \left| \sum_{p_{1}+p_{2}=p}\bar b_{p_{1}}\bar b_{p_{2}}\right|\cdot \frac{1}{(1-\tfrac{p}{M})^{s}+1},\quad \forall k\geq M\ .\\
 \end{split}
 \end{equation}
 For the remaining series in the right hand side of \eqref{eq:Fkj}, we write
  \begin{equation}
 \begin{split}
 \sum_{j=1,j\neq k}^{\infty}(\mathcal  E(j+k)+\mathcal E(|j-k|))j^{-s}&\leq  \sum_{j=1}^{\infty}\mathcal  E(j+k)j^{-s}+ \sum_{j=1}^{k-1}\mathcal  E(k-j)j^{-s}+ \sum_{j=k+1}^{\infty}\mathcal  E(j-k)j^{-s}\\
 &\leq  \sum_{j=1}^{\infty}\frac{\tilde{\mathcal  E}(j+k)}{(j+k)^{s_{\phi}}j^{s}}+ \sum_{j=1}^{k-1}\frac{\tilde{\mathcal  E}(k-j)}{(k-j)^{s_{\phi}}j^{s}}+ \sum_{j=k+1}^{\infty}\frac{\tilde{\mathcal  E}(j-k)}{(j-k)^{s_{\phi}}j^{s}}.
  \end{split}
 \end{equation}
 Since $s<s_{\phi}$, 
 \begin{equation}
 \begin{split}
 \sum_{j=1,j\neq k}^{\infty}(\mathcal  E(j+k)+\mathcal E(|j-k|))j^{-s} &\leq \frac{1}{k^{s}}\sum_{j=1}^{\infty}\frac{k^{s}\tilde{\mathcal  E}(k+j)}{(j+k)^{s}j^{s}}+ \frac{\tilde{\mathcal  E}(1)}{k^{s}}\sum_{j=1}^{k-1}\frac{k^{s}}{(k-j)^{s}j^{s}}+ \frac{1}{k^{s}}\sum_{j=k+1}^{\infty}\frac{k^{s}\tilde{\mathcal  E}(j-k)}{(j-k)^{s}j^{s}}\\
 &\leq\frac{1}{k^{s}}\left[ \sum_{j=1}^{\infty}\frac{\tilde{\mathcal  E}(k+j)}{(1+\tfrac{j}{k})^{s}j^{s}} +\tilde{\mathcal  E}(1)\gamma_{k}+\sum_{j=1}^{\infty}\frac{\tilde{\mathcal  E}(j)}{(1+\tfrac{j}{k})^{s}j^{s}}\right]
  \end{split}
 \end{equation}
 where  
 $
 \gamma_{k}
 $ is given in \eqref{eq:gammak}; see \cite{MR2718657}.
 
 Therefore, for any $k\geq M$, 
 \begin{equation}\label{eq:boundZM2}
 \begin{split}
 \sum_{j=1,j\neq k}^{\infty}(\mathcal  E(j+k)+\mathcal E(|j-k|))j^{-s} &\leq\frac{1}{k^{s}}\left[\sum_{j=1}^{\infty}\frac{ \tilde{\mathcal  E}(M+j)}{j^{s}} +\tilde{\mathcal  E}(1)\gamma_{M}+\sum_{j=1}^{\infty}\frac{\tilde{\mathcal  E}(j)}{j^{s}}\right]\\
 &\leq\frac{1}{k^{s}}\left[ \tilde{\mathcal  E}(M+1)+\frac{\tilde{\mathcal  E}(M+2)}{s-1}+\tilde{\mathcal  E}(1)\gamma_{M}+\tilde{\mathcal  E}(1)+\frac{\tilde{\mathcal  E}(2)}{s-1}\right].
   \end{split}
 \end{equation}
 Since $|c_{k,1}|_{\infty}\leq  3\sum_{j=1,j\neq k}^{\infty}|F_{k,j}|j^{-s}$, combining \eqref{eq:Fkj} with \eqref{eq:boundZM1} and \eqref{eq:boundZM2}, the thesis follows.
\eproof


\bibliographystyle{plain}

 \end{document}